\documentclass[a4paper,11pt,twoside,english]{article}
\usepackage{babel}
\usepackage{latexsym,amsfonts}
\usepackage{amsthm}
\usepackage{amsmath}
\usepackage{mathtools}
\usepackage{paralist}
\usepackage{ifthen}

\usepackage{lastpage}
\usepackage{fancyhdr}
\fancypagestyle{plain}{

  \fancyhead[L,C]{}
  \fancyhead[R]{\today}
  \fancyfoot[R,L]{}
  \fancyfoot[C]{\thepage \ of \pageref{LastPage}}
  }

\title{On convergence with respect to an ideal and a family of matrices}
\author{Jan-David Hardtke}
\date{}

\DeclareMathOperator{\co}{co}

\DeclareMathOperator{\RE}{Re}

\DeclareMathOperator{\ex}{ex}

\providecommand{\ifif}{iff }
\providecommand{\sm}{\setminus}
\providecommand{\ssq}{\subseteq}
\providecommand{\id}{\ensuremath{\mathrm{id}}}
\providecommand{\N}{\ensuremath{\mathbb{N}}}

\providecommand{\R}{\ensuremath{\mathbb{R}}}
\providecommand{\C}{\ensuremath{\mathbb{C}}}

\providecommand{\A}{\ensuremath{\mathcal{A}}}
\providecommand{\B}{\ensuremath{\mathcal{B}}}
\providecommand{\F}{\ensuremath{\mathcal{F}}}
\providecommand{\eps}{\ensuremath{\varepsilon}}
\providecommand{\cl}[2][]{\ensuremath{\overline{#2}^{#1}}}
\providecommand{\cco}[2][]{\ensuremath{\overline{\co}^{#1} #2}}

\providecommand{\keywords}[1]{
  {\let\thefootnote=\relax
  \footnote{{\em Keywords}: #1}}
  \addtocounter{footnote}{-1}
  }

\providecommand{\AMS}[1]{
  {\let\thefootnote=\relax
  \footnote{{\em AMS Subject Classification} (2010): #1}}
  \addtocounter{footnote}{-1}
  }

\providecommand{\address}{
  {\sc \noindent Department of Mathematics \\
  Freie Universit\"at Berlin \\
  Arnimallee 6, 14195 berlin \\
  Germany \\}
  }

\DeclarePairedDelimiter{\set}{\lbrace}{\rbrace}
\DeclarePairedDelimiter{\paren}{\lparen}{\rparen}

\DeclarePairedDelimiter{\abs}{\lvert}{\rvert}
\DeclarePairedDelimiter{\norm}{\lVert}{\rVert}

\theoremstyle{definition}
\newtheorem{definition}{Definition}[section]

\theoremstyle{plain}
\newtheorem{lemma}[definition]{Lemma}
\newtheorem*{lemma*}{Lemma}
\newtheorem{proposition}[definition]{Proposition}
\newtheorem{theorem}[definition]{Theorem}
\newtheorem*{theorem*}{Theorem}
\newtheorem{corollary}[definition]{Corollary}

\newenvironment{Proof}[1][\proofname]{\begin{proof}[#1] \setlength{\parindent}{0pt}}{\end{proof}}
\newenvironment{Abstract}{\centering\begin{minipage}{0.8\textwidth} \noindent \small {\sc Abstract.}}{\end{minipage}\par}

\numberwithin{equation}{section}
\usepackage{color}
\definecolor{darkgreen}{rgb}{0,0.5,0}

\newtagform{colored}[\color{blue}]{\color{blue}(}{\color{blue})}
\usetagform{colored}

\usepackage[colorlinks,linkcolor=blue,citecolor=red,urlcolor=darkgreen]{hyperref}
\providecommand{\email}{{\it E-mail address:} \href{mailto:hardtke@math.fu-berlin.de}{\tt hardtke@math.fu-berlin.de}}
\providecommand{\mr}[1]{\href{http://www.ams.org/mathscinet-getitem?mr=#1}{MR#1}}

\usepackage{amsrefs}

\begin{document}

\maketitle

\begin{Abstract}
The authors of \cite{das2} recently introduced and studied the notions of strong $A^I$-summability with respect 
to an Orlicz function $F$ and $A^I$-statistical convergence, where $A$ is a non-negative regular matrix and $I$ is
an ideal on the set of natural numbers. In this note, we will generalise these notions by replacing $A$ with a family 
of matrices and $F$ with a family of Orlicz functions or moduli and study the thus obtained convergence methods. We 
will also give an application in Banach space theory, presenting a generalisation of Simons' $\sup$-$\limsup$-theorem 
to the newly introduced convergence methods (for the case that the filter generated by the ideal $I$ has a countable 
base), continuing the work of \cite{hardtke}.
\end{Abstract}
\AMS{40C05; 40C99; 46B20}

\section{Introduction}\label{sec:intro}
Let us begin by recalling that an ideal $I$ on a non-empty set $Y$ is a non-empty set of subsets of $Y$ such that $Y\not\in I$
and $I$ is closed under the formation of subsets and finite unions. The ideal is called admissible if $\set*{y}\in I$ for each 
$y\in Y$. For example, if $Y$ is infinite then the set of all finite subsets of $Y$ forms an ideal on $Y$. If $I$ is an ideal,
then $\F(I):=\set*{Y\sm A:A\in I}$ is a filter on $Y$.\par
Now if $(x_n)_{n\in \N}$ is a sequence in a topological space $X$ and $I$ is an ideal on the set $\N$ of natural numbers then
$(x_n)_{n\in \N}$ is said to be $I$-convergent to $x\in X$ if for every neighbourhood $U$ of $x$ the set $\set*{n\in \N:x_n\not\in U}$
belongs to $I$ (equivalently, $\set*{n\in \N:x_n\in U}\in \F(I)$). In a Hausdorff space the $I$-limit is unique if it exists. It will
be denoted by $I$-$\lim x_n$. If $I_f$ is the ideal of all finite subsets of $\N$ then $I_f$-convergence is equivalent to the usual 
convergence. Thus if $I$ is admissible the usual convergence implies $I$-convergence. For a normed space $X$ the set of all $I$-convergent 
sequences in $X$ is a subspace of $X^{\N}$ and the map $(x_n) \mapsto$ $I$-$\lim x_n$ is linear. We refer the reader to \cite{kostyrko1}, 
\cite{das1}, \cite{kostyrko2} and \cite{dems} for more information on $I$-convergence.\par
Recall now that for a given infinite matrix $A=(a_{nk})_{n,k\in \N}$ with real or complex entries a sequence $s=(s_k)_{k\in \N}$
of (real or complex) numbers is said to be $A$-summable to the number $a$ provided that each of the series $\sum_{k=1}^{\infty}a_{nk}s_k$
is convergent and $\lim_{n\to \infty}\sum_{k=1}^{\infty}a_{nk}s_k=a$.\par
The matrix $A$ is called regular if every sequence that is convergent in the ordinary sense is also $A$-summable to the same limit.
A well-known theorem of Toeplitz states that $A$ is regular \ifif the following holds:
\begin{enumerate}[(i)]
\item $\sup_{n\in \N}\sum_{k=1}^{\infty}\abs{a_{nk}}<\infty$,
\item $\lim_{n\to \infty}\sum_{k=1}^{\infty}a_{nk}=1$,
\item $\lim_{n\to \infty}a_{nk}=0 \ \ \forall k\in \N$.
\end{enumerate}
Let us suppose for the moment that $A$ is regular and also non-negative (i.\,e., $a_{nk}\geq0$ for all $n,k\in \N$). We will denote by $D(s,a,\eps)$ the 
set $\set*{k\in \N: \abs{s_k-a}\geq\eps}$ for every $\eps>0$. Then $s$ is said to be $A$-statistically convergent to $a$ if for every $\eps>0$ we 
have $\lim_{n\to \infty}\sum_{k=1}^{\infty}a_{nk}\chi_{D(s,a,\eps)}(k)=0$, where the symbol $\chi_K$ denotes the characteristic function of the 
set $K\ssq \N$. If one takes $A$ to be the Ces\`aro-matrix (i.\,e., $a_{nk}=1/n$ for $k\leq n$ and $a_{nk}=0$ for $k>n$) one gets the usual notion 
of statistical convergence as it was introduced by Fast in \cite{fast}. Note that the set $I_A$ of all subsets $K\ssq \N$ for which 
$\lim_{n\to \infty}\sum_{k=1}^{\infty}a_{nk}\chi_K(k)=0$ holds, is an ideal on $\N$ and $A$-statistical convergence is nothing but convergence 
with respect to this ideal.\par
For any number $p>0$ the sequence $s$ is said to be strongly $A$-$p$-summable to $a$ provided that $\sum_{k=1}^{\infty}a_{nk}\abs{s_k-a}^p<\infty$
for all $n\in \N$ and $\lim_{n\to \infty}\sum_{k=1}^{\infty}a_{nk}\abs{s_k-a}^p=0$. The strong $A$-$p$-summability is a linear consistent summability 
method and the strong $A$-$p$-limit is uniquely determined whenever it exists. In \cite{connor1} Connor proved that $s$ is statistically convergent to 
$a$ whenever it is strongly $p$-Ces\`aro convergent to $a$ and the converse is true if $s$ is bounded. Practically the same proof as given in \cite{connor1}
still works if one replaces the Ces\`aro matrix by an arbitray non-negative regular matrix $A$. In particular, strong $A$-$p$-summability and 
$A$-statistical convergence are equivalent on bounded sequences (see also \cite{connor3}*{Theorem 8}). More information on strong matrix summability can
be found in \cite{zeller} (for the case $p=1$) or \cite{hamilton}.\par
In \cite{maddox1} Maddox proposed a generalisation of strong $A$-$p$-summability by replacing the number $p$ with a sequence $\mathbf{p}=(p_k)_{k\in \N}$
of positive numbers: the sequence s is strongly $A$-$\mathbf{p}$-summable to $a$ if $\sum_{k=1}^{\infty}a_{nk}\abs{s_k-a}^{p_k}<\infty$ for every 
$n\in \N$ and $\lim_{n\to \infty}\sum_{k=1}^{\infty}a_{nk}\abs{s_k-a}^{p_k}=0$.\par
Next, let us recall that a function $F:[0,\infty) \rightarrow [0,\infty)$ is called an Orlicz function if it is increasing, continuous, convex 
and satisfies $\lim_{t\to \infty}F(t)=\infty$ as well as $F(t)=0$ \ifif $t=0$. If we drop the convexity and replace it by the condition
$F(s+t)\leq F(s)+F(t)$ for all $s,t\geq 0$ then $F$ is called a modulus. For example, the function $F_p$ defined by $F_p(t)=t^p$ is an Orlicz
function for $p\geq 1$ and a modulus for $0<p\leq 1$. We will denote the set of all Orlicz functions by $\mathcal{O}$ and the set of all moduli
by $\mathcal{M}$.\par
Connor introduced another generalisation of strong matrix summability in \cite{connor3}: if $F$ is a modulus then $s$ is said to be strongly
$A$-summable to the limit $a$ with respect to $F$ if $\sum_{k=1}^{\infty}a_{nk}F(\abs{s_k-a})<\infty$ for all $n\in \N$ and 
$\lim_{n\to \infty}\sum_{k=1}^{\infty}a_{nk}F(\abs{s_k-a})=0$. It is shown in \cite{connor3}*{Theorem 8} that strong $A$-summability with respect 
to $F$ implies $A$-statistical convergence and that the converse holds for bounded sequences. In \cite{demirci} Demirci replaced the modulus in 
Connor's definition by an Orlicz function and studied which results carry over to this setting.\par
Another common generalised convergence method is that of almost convergence introduced by Lorentz in \cite{lorentz1}. For this we first recall that 
a Banach limit is a linear functional $L$ on the space $\ell^\infty$ of all bounded {\em real-valued} sequences such that $L$ is shift-invariant (i.\,e.,
$L((s_{n+1})_{n\in \N})=L((s_n)_{n\in \N})$), positive (i.\,e., $L((s_n)_{n\in \N})\geq0$ if $s_n\geq0$ for all $n$) and fulfils $L(1,1,\dots)=1$. The
existence of a Banach limit can be easily proved by means of the Hahn-Banach extension theorem. A sequence $s\in \ell^\infty$ is said to be almost 
convergent to $a\in \R$ if $L(s)=a$ for every Banach limit $L$.\par
It is proved in \cite{lorentz1} that almost convergence is equivalent to ``uniform Ces\`aro convergence''. More precisely, a bounded sequence 
$s=(s_k)_{k\in \N}$ in $\R$ is almost convergent to $a\in \R$ \ifif the following holds:
\begin{equation*}
 \frac{1}{n} \sum _{k=1}^{n} s_{k+i} \xrightarrow{n\to \infty} a \ \ \text{uniformly in} \ i\in \N_0,
\end{equation*}
where $\N_0=\N\cup\set*{0}$.\par
Lorentz subsequently introduced and studied the notion of $F_A$-conver\-gence by replacing the Ces\`aro-matrix with an arbitrary real-valued 
regular matrix $A$: a bounded sequence $s=(s_k)_{k\in \N}$ in $\R$ is said to be $F_A$-convergent to $a\in \R$ provided that
\begin{equation*}
\sum _{k=1}^{\infty} a_{nk}s_{k+i} \xrightarrow{n\to \infty} a \ \ \text{uniformly in} \ i\in \N_0.
\end{equation*}
Stieglitz further generalised the notion of almost convergence in the following way (cf.\,\cite{stieglitz}): consider a sequence 
$\B=(B_i)_{i\in \N_0}=((b_{nk}^{(i)})_{n,k\in \N})_{i\in \N_0}$ of matrices with entries in $\R$ or $\C$ and a bounded 
sequence $s=(s_k)_{k\in \N}$ of real or complex numbers. Then $s$ is said to be $F_{\B}$-convergent to the number $a$ 
if each of the series $\sum_{k=1}^{\infty}b_{nk}^{(i)}s_k$ with $n\in \N, i\in \N_0$ is convergent and
\begin{equation*}
\sum _{k=1}^{\infty}b_{nk}^{(i)}s_k \xrightarrow{n\to \infty} a \ \ \text{uniformly in} \ i\in \N_0.
\end{equation*}
To obtain $F_A$-convergence, take $b_{nk}^{(i)}=a_{nk-i}$ for $k>i$ and $b_{nk}^{(i)}=0$ for $k\leq i$.\par
Maddox introduced the $F_{\B}$-analogue of strong matrix summability in \cite{maddox2}. If each of the matrices $B_i$ is 
non-negative and $s=(s_k)_{k\in \N}$ is a (not necessarily bounded) sequence in $\R$ or $\C$ then $s$ is said to be strongly
$F_{\B}$-convergent to $a$ provided that
\begin{equation*}
\sum_{k=1}^{\infty}b_{nk}^{(i)}\abs{s_k-a} \xrightarrow{n\to \infty} 0 \ \ \text{uniformly in} \ i\in \N_0.
\end{equation*}
Very recently, the authors of \cite{das2} introduced the following definitions, combining matrices and ideals.
\begin{definition}[cf.\,\cite{das2}]\label{def:strong A-I stat A-I conv}
Let $A=(a_{nk})_{n,k\in \N}$ be a non-negative regular matrix, $I$ an ideal on $\N$ and $F$ an Orlicz function. Let $a$ be any 
real or complex number. A sequence $s=(s_k)_{k\in \N}$ in $\R$ or $\C$ is said to be
\begin{enumerate}[(i)]
\item strongly $A^I$-summable to $a$ with respect to $F$ if
\begin{equation*}
I\text{-}\lim\sum_{k=1}^{\infty}a_{nk}F(\abs{s_k-a})=0,
\end{equation*}
\item $A^I$-statistically convergent to $a$ if
\begin{equation*}
I\text{-}\lim\sum_{k=1}^{\infty}a_{nk}\chi_{D(s,a,\eps)}(k)=0
\end{equation*}
for every $\eps>0$.
\end{enumerate}
\end{definition}
It is proved in \cite{das2}*{Theorem 2.5} that $A^I$-summability with respect to $F$ implies $A^I$-statistical convergence (to the 
same limit) and the converse holds if the sequence $s$ is bounded and $F$ satisfies the $\Delta_2$-condition (i.\,e., there is a 
constant $K$ such that $F(2t)\leq KF(t)$ for all $t\geq 0$).\par
We would like to propose here the following three definitions that include all the above mentioned generalised convergence methods.\par
First we define a sequence $(g_n)_{n\in \N}$ of functions from a set $S$ into a generalised metric space $(X,d)$\footnote{Same as a 
metric space except that $d$ is allowed to take values in $[0,\infty]$. For example, $d(a,b)=\abs*{a-b}$ for $a,b\in [0,\infty)$,
$d(a,\infty)=d(\infty,a)=\infty$ for all $a\in [0,\infty)$ and $d(\infty,\infty)=0$ defines a generalised metric on $[0,\infty]$.}
to be uniformly convergent to
the function $g:S \rightarrow X$ along the ideal $I$ if for every $\eps>0$ there is some $E\in I$ such that for every $s\in S$
\begin{equation*}
\set*{n\in \N:d(g_n(s),g(s))\geq\eps}\ssq E
\end{equation*}
or equivalently, for every $\eps>0$ we have
\begin{equation*}
\set*{n\in \N:\sup_{i\in S}d(g_n(s),g(s))\geq\eps}\in I.
\end{equation*}
If $I=I_f$ this yields the usual definition of uniform convergence. The uniform convergence of $(g_n)_{n\in \N}$ to $g$
along $I$ clearly implies $I$-$\lim g_n(s)=g(s)$ for all $s\in S$.\par
Now for the main definition.
\begin{definition}\label{def:main def}
Let $I$ be an ideal on $\N$ and $S$ any non-empty set. Let $\B=(B_i)_{i\in S}=((b_{nk}^{(i)})_{n,k\in \N})_{i\in S}$ be
a family of (not necessarily regular) matrices with entries in $\R$ or $\C$ and $\F=(F_k^{(i)})_{k\in \N, i\in S}$ a family
in $\mathcal{M}\cup\mathcal{O}$. Suppose that there is some $i_0\in S$ such that 
\begin{equation}\label{eq:+}
\inf_{n\in \N}\sum_{k=1}^{\infty}\abs{b_{nk}^{(i_0)}}>0. \tag{+}
\end{equation}
Finally, let $s=(s_k)_{k\in \N}$ be a sequence in $\R$ or $\C$ and $a\in \R$ or $\C$.
\begin{enumerate}[(i)]
\item $s$ is said to be $\B^I$-summable to $a$ provided that each of the series $\sum_{k=1}^{\infty}b_{nk}^{(i)}s_k$ is 
convergent and 
\begin{equation*}
I\text{-}\lim\sum_{k=1}^{\infty}b_{nk}^{(i)}s_k=a \ \ \text{uniformly in} \ i\in S.
\end{equation*}
\item If each matrix $B_i$ is non-negative then $s$ is said to be strongly $\B^I$-summable to $a$ with respect to 
$\F$ if 
\begin{equation*}
I\text{-}\lim\sum_{k=1}^{\infty}b_{nk}^{(i)}F_k^{(i)}(\abs{s_k-a})=0 \ \ \text{uniformly in} \ i\in S.
\end{equation*}
\item If each $B_i$ is non-negative then $s$ is said to be $\B^I$-statistically convergent to $a$ provided that for 
every $\eps>0$ 
\begin{equation*}
I\text{-}\lim\sum_{k=1}^{\infty}b_{nk}^{(i)}\chi_{D(s,a,\eps)}(k)=0 \ \ \text{uniformly in} \ i\in S.
\end{equation*}
\end{enumerate}
\end{definition}
If $F_k^{(i)}=\id_{[0,\infty)}$ for all $k\in \N, i\in S$ in (ii) we simply speak of strong $\B^I$-summability.
Clearly, strong $\B^I$-summability to $a$ implies $\B^I$-summability to $a$ provided that $s$ is bounded,
$\sum_{k=1}^{\infty}b_{nk}^{(i)}<\infty$ for all $k\in \N, i\in S$ and 
\begin{equation*}
I\text{-}\lim\sum_{k=1}^{\infty}b_{nk}^{(i)}=1 \ \ \text{uniformly in} \ i\in S.
\end{equation*}\par
Taking $B_i=A$ and $F_k^{(i)}=F\in \mathcal{O}$ for each $i\in S$ and $k\in \N$ in (ii) and (iii) yields the definitions of
strong $A^I$-summability with respect to $F$ and of $A^I$-statistical convergence. If we take $I=I_f$ and $S=\N_0$ in (i) and 
(ii) we obtain the definitions of $F_{\B}$- and strong $F_{\B}$-convergence. Setting $I=I_f$, $B_i=A$ for 
every $i\in S$ and $F_k^{(i)}=F_{p_k}$ for all $i\in S, k\in \N$ in (ii) gives us the definition of Maddox's strong 
$A$-$\mathbf{p}$-summability.\par
Note also that if each $B_i$ is non-negative then the set $J_{\B,I}$ of all subsets $K\ssq \N$ such that
\begin{equation*}
I\text{-}\lim\sum_{k=1}^{\infty}b_{nk}^{(i)}\chi_K(k)=0 \ \ \text{uniformly in} \ i\in S,
\end{equation*}
is an ideal on $\N$ (the condition \eqref{eq:+} ensures $\N\not\in J_{\B,I}$). The $\B^I$-statistical convergence is 
nothing but the convergence with respect to $J_{\B,I}$. In the case that $B_i$ is the infinite unit matrix for each 
$i\in S$ we have $J_{\B,I}=I$.\par
In the next section we will start to investigate the above convergence methods.\par

\section{Some convergence theorems}\label{sec:conv thms}
If not otherwise stated, we will denote by $I$ an ideal on $\N$, by $\B=(B_i)_{i\in S}=((b_{nk}^{(i)})_{n,k\in \N})_{i\in S}$ a
family of real or complex matrices (where $S$ is any non-empty index set) such that there is some $i_0\in S$ with \eqref{eq:+} 
and by $\F=(F_k^{(i)})_{k\in \N, i\in S}$ a family in $\mathcal{M}\cup\mathcal{O}$. Finally, $s=(s_k)_{k\in \N}$ denotes
a sequence in and $a$ an element of $\R$ or $\C$, as in the previous section.\par
The following two propositions (wherein each $B_i$ is implicitly assumed to be non-negative) generalise the aforementioned 
results from \cite{das2}*{Theorem 2.5}. The techniques used there followed the line of \cite{connor2} while we will adopt 
the techniques from \cite{connor1}.
\begin{proposition}\label{prop:strong stat}
Suppose that $s$ is strongly $\B^I$-summable to $a$ with respect to $\F$ and that 
\begin{equation*}
L(t):=\inf\set*{F_k^{(i)}(t):k\in \N, i\in S}>0 \ \ \forall t>0.
\end{equation*}
Then $s$ is also $\B^I$-statistically convergent to $a$.
\end{proposition}

\begin{Proof}
Let $\eps, \delta>0$ be arbitrary. By assumption there is some $E\in I$ such that for all $i\in S$
\begin{equation*}
\set*{n\in \N:\sum_{k=1}^{\infty}b_{nk}^{(i)}F_k^{(i)}(\abs{s_k-a})\geq\delta L(\eps)}\ssq E.
\end{equation*}
But we have 
\begin{align*}
&\sum_{k=1}^{\infty}b_{nk}^{(i)}F_k^{(i)}(\abs{s_k-a})\geq\sum_{k=1}^{\infty}b_{nk}^{(i)}F_k^{(i)}(\abs{s_k-a})\chi_{D(s,a,\eps)}(k) \\
&\geq L(\eps)\sum_{k=1}^{\infty}b_{nk}^{(i)}\chi_{D(s,a,\eps)}(k)
\end{align*}
for all $i\in S, k\in \N$. Hence
\begin{equation*}
\set*{n\in \N:\sum_{k=1}^{\infty}b_{nk}^{(i)}\chi_{D(s,a,\eps)}(k)\geq\delta}\ssq E
\end{equation*}
for every $i\in S$ and the proof is finished.
\end{Proof}

\begin{proposition}\label{prop:stat strong}
Suppose that $s$ is bounded and $\B^I$-statistically convergent to $a$. If $\F$ is equicontinuous at $0$ and 
there exists an $A\in I$ such that
\begin{equation*}
M:=\sup\set*{\sum_{k=1}^{\infty}b_{nk}^{(i)}:n\in \N\sm A, i\in S}<\infty,
\end{equation*} 
as well as
\begin{equation*}
h(t):=\sup\set*{F_k^{(i)}(t):k\in \N, i\in S}<\infty \ \ \forall t\geq0,
\end{equation*}
then $s$ is also strongly $B^I$-summable to $a$ with respect to $\F$.
\end{proposition}

\begin{Proof}
Let $\eps>0$ be arbitray. Take $\tau>0$ with $\tau(M+h(\norm{s}_{\infty}+\abs{a}))<\eps$. Since $\F$ is equicontinuous at $0$ we can 
find a $\delta>0$ such that $F_k^{(i)}(t)\leq\tau$ for all $t\in [0,\delta]$ and all $k\in \N,i\in S$.\par
Because $s$ is $\B^I$-statistically convergent to $a$ there is some $E\in I$ such that for every $i\in S$
\begin{equation*}
\set*{n\in \N:\sum_{k=1}^{\infty}b_{nk}^{(i)}\chi_{D(s,a,\delta)}(k)\geq\tau}\ssq E.
\end{equation*}
It follows that for every $n\in \N\sm (E\cup A)$ and all $i\in S$
\begin{align*}
&\sum_{k=1}^{\infty}b_{nk}^{(i)}F_k^{(i)}(\abs{s_k-a}) \\ 
&\leq \tau\sum_{k=1}^{\infty}b_{nk}^{(i)}\chi_{\N\sm D(s,a,\delta)}(k)+\sum_{k=1}^{\infty}b_{nk}^{(i)}F_k^{(i)}(\abs{s_k-a})\chi_{D(s,a,\delta)}(k) \\
&\leq \tau M+h(\norm{s}_{\infty}+\abs{a})\sum_{k=1}^{\infty}b_{nk}^{(i)}\chi_{D(s,a,\delta)}(k)\leq\tau(M+h(\norm{s}_{\infty}+\abs{a}))<\eps
\end{align*}
and we are done.
\end{Proof}

So in particular, if $\B$ and $\F$ meet the requirements of both Proposition \ref{prop:strong stat} and Proposition \ref{prop:stat strong}
then $\B^I$-statistical convergence and strong $\B^I$-sum\-mability with respect to $\F$ coincide on bounded sequences. Note that all the 
assumptions on $\F$ are satisfied if $F_k^{(i)}=F_{p_{ki}}$ for a family $(p_{ki})_{k\in \N, i\in S}$ of positive numbers which is 
bounded and bounded away from zero.\par
If $I\ssq J_{\B,I}$, in other words, if
\begin{equation*}
I\text{-}\lim\sum_{k=1}^{\infty}b_{nk}^{(i)}\chi_A(k)=0 \ \ \mathrm{uniformly\ in} \ i\in S \ \ \forall A\in I,
\end{equation*}
then $I$-convergence implies $\B^I$-statistical convergence (to the same limit). Thus if $\B$ and $\F$ additionally satisfy the
requirements of Proposition \ref{prop:stat strong} then for bounded sequences $I$-convergence also implies strong $\B^I$-summability 
to the same limit. Concerning the consistency of ordinary $\B^I$-summability we have the following sufficient conditions which are
analogous to those of Toeplitz's theorem. We write $d_I$ for the set of all bounded sequences $(a_k)_{k\in \N}$ for which 
$\set*{k\in \N:a_k\neq 0}\in I$.
\begin{lemma}\label{lemma:BI conv cons}
Suppose that $\sum_{k=1}^{\infty}\abs{b_{nk}^{(i)}}<\infty$ for all $n\in \N, i\in S$ and 
\begin{align}
&\exists A\in I \ M:=\sup\set*{\sum_{k=1}^{\infty}\abs*{b_{nk}^{(i)}}:n\in \N\sm A, i\in S}<\infty, \label{eq:2.1} \\
&I\text{-}\lim\sum_{k=1}^{\infty}b_{nk}^{(i)}a_k=0 \ \ \mathrm{uniformly\ in} \ i\in S \ \ \forall (a_k)\in d_I, \label{eq:2.2} \\
&I\text{-}\lim\sum_{k=1}^{\infty}b_{nk}^{(i)}=1 \ \ \mathrm{uniformly\ in} \ i\in S. \label{eq:2.3}
\end{align}
Then for every bounded sequence $s=(s_n)_{n\in \N}$ in $\R$ or $\C$, if $I\text{-}\lim s_n=a$ then $s$ is also $\B^I$-summable to $a$.
\end{lemma}

\begin{Proof}
Because of \eqref{eq:2.3} we may assume $a=0$. Let $\eps>0$ be arbitrary. Since $I\text{-}\lim s_n=0$ we have 
$C:=\set*{n\in \N:\abs*{s_n}\geq\eps}\in I$ and hence by \eqref{eq:2.2} there is some $E\in I$ such that
\begin{equation*}
\set*{n\in \N:\abs*{\sum_{k=1}^{\infty}b_{nk}^{(i)}s_k\chi_C(k)}\geq\eps}\ssq E \ \ \forall i\in S.
\end{equation*}
But for all $i\in S$ and all $n\in \N\sm A$
\begin{align*}
&\abs*{\sum_{k=1}^{\infty}b_{nk}^{(i)}s_k}\leq\abs*{\sum_{k=1}^{\infty}b_{nk}^{(i)}s_k\chi_C(k)}
+\sum_{k=1}^{\infty}\abs*{b_{nk}^{(i)}}\chi_{\N\sm C}(k)\abs*{s_k} \\
&\leq\abs*{\sum_{k=1}^{\infty}b_{nk}^{(i)}s_k\chi_C(k)}+M\eps,
\end{align*}
thus
\begin{equation*}
\set*{n\in \N:\abs*{\sum_{k=1}^{\infty}b_{nk}^{(i)}s_k}\geq\eps(1+M)}\ssq E\cup A \ \ \forall i\in S
\end{equation*}
and we are done.
\end{Proof}

The next proposition is the direct generalisation of \cite{demirci2}*{Theorem 3.3} to our setting. Its proof is easy 
and moreover virtually the same as in \cite{demirci2} so it will be omitted.
\begin{proposition}\label{prop:AB-stat conv agree}
Suppose that we are given two families of non-negative matrices $\B=((b_{nk}^{(i)})_{n,k\in \N})_{i\in S}$ and 
$\A=((a_{nk}^{(i)})_{n,k\in \N})_{i\in S}$. If
\begin{equation*}
I\text{-}\lim\sum_{k=1}^{\infty}\abs*{a_{nk}^{(i)}-b_{nk}^{(i)}}=0 \ \ \mathrm{uniformly\ in} \ i\in S
\end{equation*}
then $J_{\B,I}=J_{\A,I}$.
\end{proposition}

In \cite{burgin} it was proved that a bounded (real) sequence $s$ is statistically convergent to $a$ \ifif $s$ is Ces\`aro-summable 
to $a$ and the ``variance'' $\sigma_n(s)^2:=1/n\sum_{i=1}^n\paren*{a-1/n\sum_{k=1}^ns_k}^2$ converges to $0$. The proposition below is 
a generalisation of this result. We will use the notation
\begin{equation*}
\sigma_{ni}^{\B,\F}(s):=\sum_{k=1}^{\infty}b_{nk}^{(i)}F_{ki}(\abs*{s_k-(B_i s)(n)}),
\end{equation*}
provided that each $B_i$ is non-negative.\par
First we need the following lemma, whose proof is analogous to those of Proposition \ref{prop:strong stat} and \ref{prop:stat strong}
and will therefore be omitted.
\begin{lemma}\label{lemma:aux 1}
Suppose that $\F$ and $\B$ fulfil the requirements of Proposition \ref{prop:strong stat} and Proposition \ref{prop:stat strong} and let
$y=(y_{ni})_{n\in \N,i\in S}$ be a family in $\R$ or $\C$. Put $A_{\eps,n,i}:=D(s,y_{ni},\eps)$ for all $i\in S, n\in \N$ and $\eps>0$.
Then 
\begin{equation*}
I\text{-}\lim\sum_{k=1}^{\infty}b_{nk}^{(i)}F_{ki}(\abs*{s_k-y_{ni}})=0 \ \ \mathrm{uniformly\ in} \ i\in S
\end{equation*}
implies that for every $\eps>0$
\begin{equation*}
I\text{-}\lim\sum_{k=1}^{\infty}b_{nk}^{(i)}\chi_{A_{\eps,n,i}}(k)=0 \ \ \mathrm{uniformly\ in} \ i\in S
\end{equation*}
and the converse is true if $s$ is bounded and $\sup_{i\in \N, n\in \N\sm V}\abs{y_{ni}}<\infty$ for some $V\in I$.
\end{lemma}

\begin{proposition}\label{prop:char stat conv}
Let $s$ be bounded. Under the same hypotheses as in the previous lemma and the additional assumption that
$\sum_{k=1}^{\infty}\abs{b_{nk}^{(i)}}<\infty$ for all $n\in \N, i\in S$ and
\begin{equation}\label{eq:2.4}
I\text{-}\lim\sum_{k=1}^{\infty}b_{nk}^{(i)}=1 \ \ \mathrm{uniformly\ in} \ i\in S,
\end{equation}
$s$ is $\B^I$-statistically convergent to the number $a$ \ifif $s$ is $\B^I$-summable to $a$ and $\sigma_{ni}^{\B,\F}(s)$
converges to $0$ along $I$ uniformly in $i\in S$.
\end{proposition}

\begin{Proof}
In view of Lemma \ref{lemma:aux 1} it is enough to consider the case $F_{ki}=\id_{[0,\infty)}$ for all $k\in \N, i\in S$.
We first assume that $s$ is $\B^I$-summable to $a$ and that
\begin{equation*}
I\text{-}\lim\sigma_{ni}^{\B,\F}(s)=I\text{-}\lim\sum_{k=1}^{\infty}b_{nk}^{(i)}\abs*{s_k-(B_i s)(n)}=0 \ \ \mathrm{uniformly\ in} \ i\in S.
\end{equation*}
Because of
\begin{align*}
&\sum_{k=1}^{\infty}b_{nk}^{(i)}\abs*{s_k-a}\leq\sum_{k=1}^{\infty}b_{nk}^{(i)}\abs*{s_k-(B_i s)(n)}+\sum_{k=1}^{\infty}b_{nk}^{(i)}\abs*{(B_i s)(n)-a} \\
&\leq\sigma_{ni}^{\B,\F}(s)+\abs*{(B_i s)(n)-a}M \ \ \forall n\in \N\sm A, \forall i\in S,
\end{align*}
where $A$ and $M$ are as in Proposition \ref{prop:stat strong}, it follows that $s$ is strongly $\B^I$-summable to $a$ and 
hence by Proposition \ref{prop:strong stat} it is also $\B^I$-statistically convergent to $a$.\par
Conversely, let $s$ be $\B^I$-statistically convergent to $a$. Then by Proposition \ref{prop:stat strong} $s$ 
is also strongly $\B^I$-summable to $a$ and because of our assumption \eqref{eq:2.4} it follows that $s$ is 
$\B^I$-summable to $a$. Moreover, we have
\begin{align*}
&\sigma_{ni}^{\B,\F}(s)=\sum_{k=1}^{\infty}b_{nk}^{(i)}\abs*{s_k-(B_i s)(n)}\leq\sum_{k=1}^{\infty}b_{nk}^{(i)}\abs*{s_k-a}+\sum_{k=1}^{\infty}b_{nk}^{(i)}\abs*{a-(B_i s)(n)} \\
&\leq\sum_{k=1}^{\infty}b_{nk}^{(i)}\abs*{s_k-a}+M\abs*{a-(B_i s)(n)} \ \ \forall n\in \N\sm A, \forall i\in S
\end{align*}
and hence $\sigma_{ni}^{\B,\F}(s)$ converges to $0$ along $I$ uniformly in $i\in S$.
\end{Proof}

According to \cite{lorentz1}*{Theorem 2}, for any regular matrix $A$ the $F_A$-covergence of a sequence implies its 
almost convergence to the same limit and by \cite{lorentz1}*{Theorem 3} the converse is true if $A$ satisfies 
$\lim_{n\to \infty}\sum_{k=1}^{\infty}\abs*{a_{nk}-a_{nk+1}}=0$. The following two results are generalisations of
these facts. Their proofs remain virtually the same and will not be given here.
\begin{proposition}\label{prop:strong AI conv almost conv}
Let $A=(a_{nk})_{n,k\in \N}$ be an infinite matrix in $\R$ such that $\sum_{k=1}^{\infty}\abs*{a_{nk}}<\infty$ for
all $n\in \N$ and $I\text{-}\lim\sum_{k=1}^{\infty}a_{nk}=1$. Put $\A=((a_{nk}^{(i)})_{n,k\in \N})_{i\in \N_0}$, where 
$a_{nk}^{(i)}=a_{nk-i}$ for $k>i$ and $a_{nk}^{(i)}=0$ for $k\leq i$.\par
Let $s\in \ell^{\infty}$ be $\A^I$-summable to the value $a$. Then $s$ is also almost convergent to $a$.
\end{proposition}

\begin{theorem}
Let $A$ and $\A$ be as in the previous proposition but assume additionally that $I\text{-}\lim a_{nk}=0$ for every $k\in \N$,
$\sup_{n\in \N\sm V}\sum_{k=1}^{\infty}\abs*{a_{nk}}<\infty$ for some $V\in I$ and
\begin{equation*}
I\text{-}\lim\sum_{k=1}^{\infty}\abs*{a_{nk}-a_{nk+1}}=0.
\end{equation*}
Let $C$ be the Ces\`aro-matrix and suppose that the family $\mathcal{C}$ arises from $C$ as $\A$ from $A$. Suppose further that 
the ideal $I$ is admissible and that $J$ is another ideal. Let $s\in \ell^{\infty}$ be $\mathcal{C}^J$-summable to the value $a$. 
Then $s$ is also $\A^I$-summable to $a$.
\end{theorem}

In \cite{dems} the notion of $I$-Cauchy sequences in arbitrary metric spaces, which generalises the notion of statistically Cauchy 
sequences of Fridy (cf.\,\cite{fridy1}), was introduced. A sequence $(x_n)_{n\in \N}$ in a metric space $(X,d)$ is said to be an $I$-Cauchy 
sequence if for every $\eps>0$ there is some $k\in \N$ such that $\set*{n\in \N: d(x_n,x_k)\geq\eps}\in I$. For $I=I_f$ this yields
just an equivalent formulation of the notion of an ordinary Cauchy sequence. Fridy's notion of statistically Cauchy sequences is 
obtained by taking $I=J_{C,I_f}$, where $C$ is the Ces\`aro-matrix. It was proved in \cite{dems} that every $I$-convergent sequence
is $I$-Cauchy (cf.\,\cite{dems}*{Proposition 1}) and that, in the case of an admissible ideal $I$, the metric space $(X,d)$ is 
complete \ifif every $I$-Cauchy sequence in $(X,d)$ is $I$-convergent (cf.\,\cite{dems}*{Theorem 2}). The proof of 
\cite{dems}*{Theorem 2} also shows that every $I$-convergent sequence possesses a subsequence which is convergent in the 
ordinary sense.\par
In \cite{fridy1} it was proved that a sequence of numbers is statistically convergent \ifif it is statistically
Cauchy, but a third equivalent condition was obtained there as well, namely a number sequence $(s_n)_{n\in \N}$ is statistically
convergent \ifif there is a sequence $(t_n)_{n\in \N}$ which is convergent in the usual sense and coincides ``almost everywhere''
with $(s_n)_{n\in \N}$, which in our notation means precisely $\set*{n\in \N:s_n\neq t_n}\in J_{C,I_f}$.\par
It is clear that for any two sequences $(x_n)_{n\in \N}, (y_n)_{n\in \N}$ in an arbitrary topological space, if $(y_n)_{n\in \N}$
is $I$-convergent and $\set*{n\in \N:x_n\neq y_n}\in I$ then $(x_n)_{n\in \N}$ is also $I$-convergent. For the case of $\B^I$-statistical
convergence of sequences of numbers we can prove a converse result provided that $\F(I)$ has a countable base that fulfils a certain 
condition with respect to the matrix-family $\B$. The proof uses the basic ideas from \cite{fridy1}.

\begin{theorem}\label{thm:char BI-stat conv}
Let $I$ be an admissible ideal with $I\ssq J_{\B,I}$ such that there is an increasing sequence $(B_m)_{m\in \N}$ in $I$ 
for which $\set*{\N\sm B_m:m\in \N}$ forms a base of $\F(I)$ and
\begin{equation}\label{eq:2.5}
\sup_{i\in S}\sup_{n\in B_m}\sum_{k=1}^{\infty}b_{nk}^{(i)}\chi_{\N\sm B_m}(k) \xrightarrow{m\to \infty} 0.
\end{equation}
Then the sequence $s=(s_n)_{n\in \N}$ is $\B^I$-statistically convergent to $a$ \ifif there is a sequence $(t_n)_{n\in \N}$
which is $I$-convergent to $a$ and fulfils $\set*{n\in \N:s_n\neq t_n}\in J_{\B,I}$.
\end{theorem}

\begin{Proof}
We only have to show the necessity. So let $s$ be $\B^I$-statistically convergent to $a$. Put $\eps_m=2^{-m}$ and 
$A_m=\set*{k\in \N:\abs*{s_k-a}\geq\eps_m}$ for every $m\in \N$. Then for every $m\in \N$ there exists a set $E_m\in I$ 
such that
\begin{equation}\label{eq:2.6}
\set*{n\in \N:\sum_{k=1}^{\infty}b_{nk}^{(i)}\chi_{A_m}(k)\geq\eps_m}\ssq E_m \ \ \forall i\in S
\end{equation}
and by \eqref{eq:2.5} we can find a strictly increasing sequence $(M_p)_{p\in \N}$ in $\N$ such that 
\begin{equation}\label{eq:2.7}
\sup_{i\in S}\sup_{n\in B_{M_p}}\sum_{k=1}^{\infty}b_{nk}^{(i)}\chi_{\N\sm B_{M_p}}(k)\leq\eps_p \ \ \forall p\in \N.
\end{equation}
Next we fix a strictly increasing sequence $(p_m)_{m\in \N}$ in $\N$ such that $E_m\ssq B_{M_{p_m}}$ for every $m\in \N$. We write
$F_m$ for $B_{M_{p_m}}$. Then $F_m\ssq F_{m+1}$ and $\bigcup_{m=1}^{\infty}F_m=\N$.\par
Let $m(k)=\min\set*{m\in \N: k\in F_m}$ for every $k\in \N$ and put
\begin{equation*}
t_k=
\begin{cases}
s_k \ \mathrm{if} \ k\not\in A_{m(k)} \\
a \ \mathrm{if} \ k\in A_{m(k)}.
\end{cases}
\end{equation*}
It is easily checked that $\set*{k\in \N:\abs*{t_k-a}\geq\eps_m}\ssq F_m$ for every $m$ and hence $(t_k)_{k\in \N}$ 
is $I$-convergent to $a$.\par
Now it remains to show $C:=\set*{k\in \N:s_k\neq t_k}\in J_{\B,I}$. To this end, fix $\eps>0$ and choose $m$ such that 
$\sum_{l=m+1}^{\infty}\eps_l\leq\eps/3$ and $\eps_{p_m}\leq\eps/3$.\par
Since $I\ssq J_{\B,I}$ we can find $E\in I$ with
\begin{equation}\label{eq:2.8}
\set*{n\in \N:\sum_{k=1}^{\infty}b_{nk}^{(i)}\chi_{F_m}(k)\geq\frac{\eps}{3}}\ssq E \ \ \forall i\in S.
\end{equation}
Then $F_m\cup E\in I$ and for every $n\in \N\sm (F_m\cup E)$ and each $i\in S$ we have $m(n)>m$ and
\begin{align*}
&\sum_{k=1}^{\infty}b_{nk}^{(i)}\chi_C(k)=\sum_{k=1}^{\infty}b_{nk}^{(i)}\chi_{C\cap F_m}(k)+\sum_{k=1}^{\infty}b_{nk}^{(i)}\chi_{C\cap (\N\sm F_m)}(k) \\
&\stackrel{\eqref{eq:2.8}}{<}\frac{\eps}{3}+\sum_{k=1}^{\infty}b_{nk}^{(i)}\chi_{C\cap (\N\sm F_{m(n)})}(k)
+\sum_{k=1}^{\infty}b_{nk}^{(i)}\chi_{C\cap (F_{m(n)}\sm F_m)}(k) \\
&\stackrel{\eqref{eq:2.7}}{\leq}\frac{\eps}{3}+\eps_{p_{m(n)}}+\sum_{l=m+1}^{m(n)}\sum_{k=1}^{\infty}b_{nk}^{(i)}\chi_{C\cap(F_l\sm F_{l-1})}(k) \\
&\leq\frac{\eps}{3}+\eps_{p_m}+\sum_{l=m+1}^{m(n)}\sum_{k=1}^{\infty}b_{nk}^{(i)}\chi_{A_l}(k) \\
&\stackrel{\eqref{eq:2.6}}{\leq}\frac{2}{3}\eps+\sum_{l=m+1}^{m(n)}\eps_l\leq\eps,
\end{align*}
which completes the proof.
\end{Proof}
Note that condition \eqref{eq:2.5} is in particular satisfied for $B_m=\set*{1,\dots,m}$ if $I=I_f$ and each 
$B_i$ is a lower triangular matrix.

Making use of his aforementioned characterisation of statistical convergence, Fridy further proved in \cite{fridy1} the following 
Tauberian theorem for statistical convergence: a statistically convergent sequence $(s_n)_{n\in \N}$ which satisfies 
$\abs*{s_n-s_{n+1}}=O(1/n)$ for $n\to \infty$ is convergent in the ordinary sense. It is not too difficult to obtain the following
slightly more general result by modifying the proof from \cite{fridy1} accordingly (there the functions $\varphi, \psi$ and $h$
below are simply $\varphi(x)=1/x=\psi(x)$ and $h(x)=x(1+x)^{-1}$). For the sake of brevity, we skip the details.
\begin{theorem}\label{thm:tauberian cond}
Let $I$ be an admissible ideal and $A=(a_{nk})_{n,k\geq 1}$ a lower triangular matrix such that $I\text{-}\lim\sum_{k=1}^na_{nk}=1$
and $I\text{-}\lim a_{nk}=0$ for every $k\in \N$. Suppose that $\varphi,\psi$ and $h$ are functions from $[0,\infty)$ into itself
such that $\varphi$ is decreasing on $(0,\infty)$, $\min_{k=1,\dots,n}a_{nk}\geq\psi(n)$ for every $n\in \N$, $I\text{-}\lim x_n=0$
whenever $I\text{-}\lim h(x_n)=0$, and 
\begin{equation*}
x\psi(x+y)\geq h(x\varphi(y)) \ \ \forall x,y\geq 0.
\end{equation*}
Let $(s_n)_{n\in \N}$ and $(t_n)_{n\in \N}$ be number sequences such that $\lim_{n\to \infty} t_n=0$, $\set*{n\in \N:s_n\neq t_n}\in J_{A,I}$
and $\abs*{s_n-s_{n+1}}=O(\varphi(n))$ for $n\to \infty$. Then $I\text{-}\lim s_n=0$.
\end{theorem}
Combining the Theorems \ref{thm:char BI-stat conv} and \ref{thm:tauberian cond} we get the following corollary.
\begin{corollary}\label{cor:tauberian cond}
Under the same general hypothesis as in Theorem \ref{thm:tauberian cond} with $I=I_f$, if $(s_n)_{n\in \N}$ is a
sequence which is $A$-statistically convergent to the number $a$ and fulfils $\abs*{s_n-s_{n+1}}=O(\varphi(n))$ 
for $n\to \infty$, then $(s_n)_{n\in \N}$ is convergent to $a$ in the usual sense.
\end{corollary}

\section{Limit superior and limit inferior}\label{sec:limsup}
In \cite{demirci0} Demirci introduced the concepts of limit superior and limit inferior with respect to an ideal $I$
on $\N$, generalising the notions of statistical limit superior and limit inferior from \cite{fridy2}. For a sequence 
$(s_n)_{n\in \N}$ in $\R$ put
\begin{align*}
&I\text{-}\limsup s_n:=\sup\set*{t\in \R:\set*{n\in \N:s_n>t}\not\in I}, \\
&I\text{-}\liminf s_n:=\inf\set*{t\in \R:\set*{n\in \N:s_n<t}\not\in I}.
\end{align*}
The same definitions were independently introduced by the authors of \cite{kostyrko2}. Note that since $(s_n)_{n\in \N}$ is 
not assumed to be bounded, it can happen that these values are $\infty$ or $-\infty$. If $I=I_f$ the above definitions are 
equivalent to the usual definitions of limit superior and limit inferior. It is proved in \cite{demirci0} (and in \cite{kostyrko2} 
as well) that $I\text{-}\liminf s_n\leq I\text{-}\limsup s_n$ and that $(s_n)_{n\in \N}$ is $I$-convergent to 
$a\in \R$ \ifif $I\text{-}\liminf s_n=a=I\text{-}\limsup s_n$ (cf.\,\cite{demirci0}*{Theorems 3 and 4} or 
\cite{kostyrko2}*{Theorems 3.2 and 3.4}).\par
Let us also remark that 
\begin{equation*}
I\text{-}\limsup s_n=\inf_{A\in I}\sup\set*{s_n:n\in \N\sm A}
\end{equation*}
and
\begin{equation*}
I\text{-}\liminf s_n=\sup_{A\in I}\inf\set*{s_n:n\in \N\sm A},
\end{equation*}
as is easily checked.\par
In \cite{fridy2}*{Lemma on p.3628} necessary and sufficient conditions for a real matrix $A$ to satisfy the inequality
$\limsup Ax\leq \operatorname{st-lim\,sup}x$ for all $x\in \ell^{\infty}$ were obtained (here, $\operatorname{st-lim\,sup}x$ 
denotes the aforementioned statistical limit superior that was introduced in \cite{fridy2}, in our terminology it is nothing 
but the limit superior with respect to the ideal $J_{C,I_f}$, where $C$ is the Ces\`aro-matrix).\par
Later, Demirci gave a more general necessity result concerning the $I$-limit superior and the $I$-limit inferior
(cf.\,\cite{demirci0}*{Corollary 1}). The following proposition is a further generalisation of this result while its proof 
follows the lines from \cite{fridy2}.
\begin{proposition}\label{prop:A limsup}
Let $I,J$ be ideals on $\N$ and $A=(a_{nk})_{n,k\in \N}$ an infinite matrix in $\R$ such that the following conditions are satisfied:
\begin{align}
&\sum_{k=1}^{\infty}\abs{a_{nk}}<\infty \ \ \forall n\in \N, \label{eq:3.1} \\
&I\text{-}\lim\sum_{k=1}^{\infty}\abs{a_{nk}}=1=I\text{-}\lim\sum_{k=1}^{\infty}a_{nk}, \label{eq:3.2} \\
&I\text{-}\lim\sum_{k=1}^{\infty}\abs{a_{nk}}\chi_E(k)=0 \ \ \forall E\in J. \label{eq:3.3}
\end{align}
Then 
\begin{equation*}
I\text{-}\limsup As\leq J\text{-}\limsup s \ \ \forall s\in \ell^{\infty}
\end{equation*}
as well as
\begin{equation*}
I\text{-}\liminf As\geq J\text{-}\liminf s \ \ \forall s\in \ell^{\infty}.
\end{equation*}
\end{proposition}

\begin{Proof}
Let $s=(s_n)_{n\in \N}\in \ell^{\infty}$ be arbitrary and put $b=J\text{-}\limsup s$. Since $s$ is bounded we have $b\in \R$.
Also, fix an arbitrary $\eps>0$. Then by \cite{demirci0}*{Theorem 1} (or \cite{kostyrko2}*{Theorem 3.1}) we have 
$E:=\set*{n\in \N:s_n>b+\eps}\in J$. We put $F=\N\sm E$.\par
For every $a\in \R$ set $a^+=\max\set*{a,0}$ and $a^-=\max\set*{-a,0}$, as in \cite{fridy2}. Note that $a=a^+-a^+$ and 
$\abs{a}=a^++a^-$.\par
Then for every $n\in \N$
\begin{align*}
&(As)(n)=\sum_{k=1}^{\infty}a_{nk}s_k=\sum_{k=1}^{\infty}a_{nk}^+\chi_E(k)s_k+\sum_{k=1}^{\infty}a_{nk}^+\chi_F(k)s_k-\sum_{k=1}^{\infty}a_{nk}^-s_k \\
&\leq \norm{s}_{\infty}\sum_{k=1}^{\infty}\abs{a_{nk}}\chi_E(k)+(b+\eps)\sum_{k=1}^{\infty}a_{nk}^+\chi_F(k)+\frac{1}{2}\norm{s}_{\infty}\sum_{k=1}^{\infty}(\abs{a_{nk}}-a_{nk}) \\
&=\norm{s}_{\infty}\sum_{k=1}^{\infty}\abs{a_{nk}}\chi_E(k)+\frac{1}{2}\norm{s}_{\infty}\sum_{k=1}^{\infty}(\abs{a_{nk}}-a_{nk}) \\
&+\frac{b+\eps}{2}\paren*{\sum_{k=1}^{\infty}(\abs{a_{nk}}+a_{nk})(1-\chi_E(k))}.
\end{align*}
Because of $E\in J$ and the assumptions \eqref{eq:3.2} and \eqref{eq:3.3} the $I$-limit of the right-hand side of the above inequality 
is equal to $b+\eps$. Together with the obvious monotonicity of $I$-$\limsup$ it follows that $I\text{-}\limsup As\leq b+\eps$. 
Since $\eps>0$ was arbitrary, the proof is finished.\par
The second statement follows from the first one by multiplication with $-1$.
\end{Proof}

It was also proved in \cite{fridy2} that a sequence of real numbers which is bounded above and Ces\`aro-summable to its
statistical limit superior is statistically convergent (cf.\,\cite{fridy2}*{Theorem 5}). It is possible to modify the proof
of \cite{fridy2} to obtain the following more general result. We use the same notation as in the previous section.
\begin{theorem}\label{thm:limsup stat}
Suppose that each $B_i$ is non-negative, $\sum_{k=1}^{\infty}b_{nk}^{(i)}<\infty$ for all $n\in \N, i\in S$ and
\begin{equation}\label{eq:3.4}
I\text{-}\lim\sum_{k=1}^{\infty}b_{nk}^{(i)}=1 \ \ \mathrm{uniformly \  in} \ i\in S.
\end{equation}
If $s=(s_n)_{n\in \N}$ is a bounded sequence of real numbers and $a\in \R$ such that $s$ is $\B^I$-summable to $a$
and $J_{\B,I}$-$\limsup s=a$ or $J_{\B,I}$-$\liminf s=a$ then $s$ is $\B^I$-statistically convergent to $a$.
\end{theorem}

\begin{Proof}
It is enough to prove the statement for the case $J_{\B,I}$-$\limsup s=a$. Suppose that $s$ is not $B^I$-statistically convergent 
to $a$. Then $J_{\B,I}$-$\liminf s<a$ and hence there must be some $t<a$ such that $E:=\set*{n\in \N:s_n<t}\not\in J_{\B,I}$.
Consequently, there exists a $d>0$ such that 
\begin{equation}\label{eq:3.5}
A:=\set*{n\in \N:\sup_{i\in S}\sum_{k=1}^{\infty}b_{nk}^{(i)}\chi_E(k)\geq d}\not\in I.
\end{equation}
Fix an arbitrary $\eps>0$ and put $F:=\set*{n\in \N:t\leq s_n\leq a+\eps}$ and $G:=\set*{n\in \N:s_n>a+\eps}$. Take $\delta\in (0,\eps)$
with $\delta\abs{a+\eps}\leq\eps$. By our assumption \eqref{eq:3.4} we have
\begin{equation*}
C:=\set*{n\in \N:\sup_{i\in S}\abs*{\sum_{k=1}^{\infty}b_{nk}^{(i)}-1}\geq\delta}\in I.
\end{equation*}
It follows from \cite{demirci0}*{Theorem 1} that $G\in J_{\B,I}$ and hence
\begin{equation*}
D:=\set*{n\in \N:\sup_{i\in S}\sum_{k=1}^{\infty}b_{nk}^{(i)}\chi_G(k)\geq\delta}\in I.
\end{equation*}
Now let $n\in H:=A\cap(\N\sm(C\cup D))$ be arbitrary. Since $n\in A$ there is some $i\in S$ such that 
$\sum_{k=1}^{\infty}b_{nk}^{(i)}\chi_E(k)>d/2$. Write $M=\norm{s}_{\infty}$. It then follows from the definitions
of the sets $E,F,G,C$ and $D$ and the choice of $\delta$ that
\begin{align*}
&\sum_{k=1}^{\infty}b_{nk}^{(i)}s_k=\sum_{k=1}^{\infty}b_{nk}^{(i)}s_k\chi_E(k)+\sum_{k=1}^{\infty}b_{nk}^{(i)}s_k\chi_F(k)
+\sum_{k=1}^{\infty}b_{nk}^{(i)}s_k\chi_G(k) \\
&\leq t\sum_{k=1}^{\infty}b_{nk}^{(i)}\chi_E(k)+(a+\eps)\sum_{k=1}^{\infty}b_{nk}^{(i)}\chi_F(k)+M\delta \\
&=t\sum_{k=1}^{\infty}b_{nk}^{(i)}\chi_E(k)+M\delta+(a+\eps)\sum_{k=1}^{\infty}b_{nk}^{(i)}(1-\chi_E(k)-\chi_G(k)) \\
&\leq t\sum_{k=1}^{\infty}b_{nk}^{(i)}\chi_E(k)+M\delta+(a+\eps)\paren*{1-\sum_{k=1}^{\infty}b_{nk}^{(i)}\chi_E(k)} \\
&+\abs*{a+\eps}\paren*{\sum_{k=1}^{\infty}b_{nk}^{(i)}\chi_G(k)+\abs*{\sum_{k=1}^{\infty}b_{nk}^{(i)}-1}} \\
&\leq a+\eps+M\eps+(t-a-\eps)\sum_{k=1}^{\infty}b_{nk}^{(i)}\chi_E(k)+2\abs*{a+\eps}\delta \\
&<a+\eps(M+3)+(t-a-\eps)\frac{d}{2}.
\end{align*}
Thus we have 
\begin{equation*}
\sup_{i\in S}\abs*{a-\sum_{k=1}^{\infty}b_{nk}^{(i)}s_k}>\frac{d}{2}(a+\eps-t)-\eps(M+3) \ \ \forall n\in H.
\end{equation*}
Suppose that
\begin{equation*}
h:=I\text{-}\limsup \sup_{i\in S}\abs*{a-\sum_{k=1}^{\infty}b_{nk}^{(i)}s_k}<\frac{d}{2}(a+\eps-t)-\eps(M+3).
\end{equation*}
Then it would follow that $H\in I$. But $C,D\in I$ and hence 
\begin{equation*}
A=H\cup(C\cap A)\cup(D\cap A)\in I,
\end{equation*}
contradicting \eqref{eq:3.5}.\par
Thus $h\geq\frac{d}{2}(a+\eps-t)-\eps(M+3)$ and since $\eps>0$ was arbitrary we get $h\geq (a-t)d/2>0$ and hence 
$s$ is not $B^I$-summable to $a$.
\end{Proof}

We conclude this section with a lemma that will be needed later and may also be of independent interest. First we need 
one more definition: a number sequence $(s_n)_{n\in \N}$ is called $I$-bounded if there is a constant $K>0$ such that 
$\set*{n\in \N:\abs*{s_n}>K}\in I$. Note that $I$-convergent sequences are $I$-bounded and that the $I$-boundedness of 
$(s_n)_{n\in \N}$ implies that $I\text{-}\limsup s_n$ and $I\text{-}\liminf s_n$ are finite.
\begin{lemma}\label{lemma:I limsup}
For any ideal $I$ on $\N$ and all $I$-bounded sequences $(s_n)_{n\in \N}$ and $(t_n)_{n\in \N}$ in $\R$ the inequalities
\begin{align*}
&I\text{-}\limsup (s_n+t_n)\leq I\text{-}\limsup s_n+I\text{-}\limsup t_n \ \ \mathrm{and} \\
&I\text{-}\liminf (s_n+t_n)\geq I\text{-}\liminf s_n+I\text{-}\liminf t_n
\end{align*}
hold. If one of the sequences is $I$-convergent, then equality holds.
\end{lemma}

\begin{Proof}
It is enough to prove the statement for the $I\text{-}\limsup$. Let $a=I\text{-}\limsup s_n$ and $b=I\text{-}\limsup t_n$.
If $u,v\in \R$ such that $u>a$ and $v>b$ then $A:=\set*{n\in \N:s_n>u}\in I$ and $B:=\set*{n\in \N:t_n>v}\in I$. 
Hence $A\cup B\in I$. But
\begin{equation*}
C:=\set*{n\in \N:s_n+t_n>u+v}\ssq A\cup B,
\end{equation*}
thus $C\in I$.\par
If $I\text{-}\limsup(s_n+t_n)>u+v$ then there would be some $\eta>u+v$ such that $\set*{n\in \N:s_n+t_n>\eta}\not\in I$,
which would imply $C\not\in I$. Thus we must have $I\text{-}\limsup(s_n+t_n)\leq u+v$. Since $u>a$ and $v>b$ were arbitrary
it follows that $I\text{-}\limsup(s_n+t_n)\leq a+b$.\par
Now suppose that $(s_n)_{n\in \N}$ is $I$-convergent to $a$ and fix an arbitrary $\eps>0$. Put $D:=\set*{n\in \N:s_n+t_n>a+b-\eps}$,
$E:=\set*{n\in \N:s_n>a-\eps/2}$ and $F:=\set*{n\in \N:t_n>b-\eps/2}$.\par
By \cite{demirci0}*{Theorem 1} $F\not\in I$ and because of $I\text{-}\lim s_n=a$ we have $\N\sm E\in I$, i.\,e., $E\in \F(I)$.\par
If $E\cap F\in I$ then $(\N\sm E)\cup (\N\sm F)\in \F(I)$ and hence $(\N\sm F)\cap E=((\N\sm E)\cup (\N\sm F))\cap E\in \F(I)$, thus
$\N\sm F\in \F(I)$, contradicting the fact that $F\not\in I$.\par
So we must have $E\cap F\not\in I$ and since $E\cap F\ssq D$ it follows that $D\not\in I$, which implies 
$I\text{-}\limsup(s_n+t_n)\geq a+b-\eps$. Letting $\eps\to 0$ completes the proof.
\end{Proof}

\section{Cluster points}
Fridy (\cite{fridy}) defined and studied statistical cluster points and statistical limit points of a sequence. These 
concepts were later generalised by the authors of \cite{kostyrko1} to an arbitrary admissible ideal $I$. Consider a 
sequence $(x_n)_{n\in \N}$ in a metric space $(X,d)$. An element $x\in X$ is called an $I$-cluster point of 
$(x_n)_{n\in \N}$ if $\set*{n\in \N: d(x_n,x)<\eps}\not\in I$ for every $\eps>0$ and it is called an $I$-limit point 
of $(x_n)_{n\in \N}$ if there is a subsequence $(x_{n_k})_{k\in \N}$ with $\set*{n_k:k\in \N}\not\in I$ that converges 
to $x$. For $I=I_f$, both notions are equivalent to the usual notion of cluster points. Every $I$-limit point is also
an $I$-cluster point of $(x_n)_{n\in \N}$ (cf.\,\cite{kostyrko1}*{Proposition 4.1}) but the converse is not true in general.
It was shown in \cite{kostyrko2}*{Theorem 3.5} that a bounded sequence $(s_n)_{n\in \N}$ in $\R$ always possesses an 
$I$-cluster point and that the $I\text{-}\limsup$ and the $I\text{-}\liminf$ of the sequence is the greatest respectively 
the smallest of them. It is easily observed that the same proof still works if the sequence is only $I$-bounded.\par
Concerning $J_{\B,I}$-cluster points, we can give the following characterisation.
\begin{proposition}\label{prop:char cluster}
Suppose that $\sup_{n\in \N,i\in S}\sum_{k=1}^{\infty}b_{nk}^{(i)}<\infty$ and
\begin{equation}\label{eq:4.1}
I\text{-}\lim\sum_{k=1}^{\infty}b_{nk}^{(i)}=1 \ \ \mathrm{uniformly\ in} \ i\in S.
\end{equation}
Then $a$ is a $J_{\B,I}$-cluster point of $s=(s_n)_{n\in \N}$ \ifif for every $\eps>0$ 
\begin{equation*}
I\text{-}\liminf\inf_{i\in S}\sum_{k=1}^{\infty}b_{nk}^{(i)}\chi_{D(s,a,\eps)}(k)<1.
\end{equation*}
\end{proposition}

\begin{Proof}
Put $A_{\eps}=D(s,a,\eps)$ and $B_{\eps}=\N\sm A_{\eps}$ for every $\eps>0$. By definition, $a$ is a $J_{\B,I}$-cluster 
point of $s$ \ifif $B_{\eps}\not\in J_{\B,I}$ for every $\eps>0$ which is the case \ifif
\begin{equation*}
I\text{-}\limsup\sup_{i\in S}\sum_{k=1}^{\infty}b_{nk}^{(i)}\chi_{B_{\eps}}(k)>0.
\end{equation*}
But $\sum_{k=1}^{\infty}b_{nk}^{(i)}\chi_{B_{\eps}}(k)=\sum_{k=1}^{\infty}b_{nk}^{(i)}
-\sum_{k=1}^{\infty}b_{nk}^{(i)}\chi_{A_{\eps}}(k)$, so because of \eqref{eq:4.1} and Lemma \ref{lemma:I limsup} it follows
that $a$ is a $J_{\B,I}$-cluster point of $s$ \ifif 
\begin{align*}
&I\text{-}\limsup\sup_{i\in S}\paren*{1-\sum_{k=1}^{\infty}b_{nk}^{(i)}\chi_{A_{\eps}}(k)}>0 \ \iff \\
&1-I\text{-}\liminf\inf_{i\in S}\sum_{k=1}^{\infty}b_{nk}^{(i)}\chi_{A_{\eps}}(k)>0
\end{align*}
and the proof is finished.
\end{Proof}

This characterisation yields the following sufficient condition for a $J_{\B,I}$-cluster point.
\begin{corollary}\label{cor:cluster}
Under the same assumptions as in the previous proposition, if $\F=(F_k^{(i)})_{k\in \N, i\in S}$ is a family in 
$\mathcal{M}\cup\mathcal{O}$ such that
\begin{align*}
&L(t):=\inf\set*{F_k^{(i)}(t):k\in \N, i\in S}>0 \ \ \forall t>0 \ \ \mathrm{and} \\
&I\text{-}\liminf\inf_{i\in S}\sum_{k=1}^{\infty}b_{nk}^{(i)}F_k^{(i)}(\abs*{s_k-a})=0,
\end{align*}
then $a$ is a $J_{\B,I}$-cluster point of $s$.
\end{corollary}

\begin{Proof}
For every $\eps>0$ and all $i\in S, n\in \N$ we have
\begin{align*}
&\sum_{k=1}^{\infty}b_{nk}^{(i)}F_k^{(i)}(\abs*{s_k-a})\geq\sum_{k=1}^{\infty}b_{nk}^{(i)}F_k^{(i)}(\abs*{s_k-a})\chi_{D(s,a,\eps)}(k) \\
&\geq L(\eps)\sum_{k=1}^{\infty}b_{nk}^{(i)}\chi_{D(s,a,\eps)}(k)
\end{align*}
and thus it follows from the assumptions that
\begin{equation*}
I\text{-}\liminf\inf_{i\in S}\sum_{k=1}^{\infty}b_{nk}^{(i)}\chi_{D(s,a,\eps)}(k)=0<1 \ \ \forall \eps>0.
\end{equation*}
Hence by the previous proposition, $a$ is a $J_{\B,I}$-cluster point of $s$.
\end{Proof}

\section{Pre-Cauchy sequences}\label{sec:pre cauchy}
The authors of \cite{connor4} introduced the notion of statistically pre-Cauchy sequences. The sequence $s=(s_k)_{k\in \N}$
is called a statistically pre-Cauchy sequence if $\lim_{n\to \infty}1/n^2\abs*{\set*{(i,j)\in \set*{1,\dots,n}^2:\abs{s_i-s_j}\geq\eps}}=0$
for every $\eps>0$. They show that a statistically convergent sequence is statistically pre-Cauchy and that the converse is
not true in general but under certain additional assumptions. It is further proved that $s$ is statistically pre-Cauchy if
\begin{equation*}
\lim_{n\to \infty}\frac{1}{n^2}\sum_{i=1}^n\sum_{j=1}^n\abs*{s_i-s_j}=0
\end{equation*}
and that the converse is true if $s$ is bounded (cf.\,\cite{connor4}*{Theorem 3}).\par
We propose the following generalisation of the definition of statistically pre-Cauchy sequences to our setting.
\begin{definition}\label{def:pre cauchy}
If each $B_i$ is non-negative, a sequence $s=(s_k)_{k\in \N}$ of real or complex numbers is called a $\B^I$-statistically
pre-Cauchy sequence if for every $\eps>0$
\begin{equation*}
I\text{-}\lim\sum_{k=1}^{\infty}\sum_{l=1}^{\infty}b_{nk}^{(i)}b_{nl}^{(i)}\chi_{D(s,\eps)}(k,l)=0 \ \ \mathrm{uniformly\ in} \ i\in S,
\end{equation*}
where $D(s,\eps):=\set*{(k,l)\in \N^2:\abs*{s_k-s_l}\geq\eps}$.
\end{definition}

First we show that, under an additional assumption on $\B$, $\B^I$-statistically convergent sequences are 
$\B^I$-statistically pre-Cauchy.
\begin{lemma}\label{lemma:stat conv pre cauchy}
Suppose that $s$ is $\B^I$-statistically convergent and 
\begin{equation*}
\exists A\in I\ M:=\sup\set*{\sum_{k=1}^{\infty}b_{nk}^{(i)}:n\in \N\sm A, i\in S}<\infty.
\end{equation*}
Then $s$ is a $\B^I$-statistically pre-Cauchy sequence.
\end{lemma}

\begin{Proof}
Say $s$ is $\B^I$-statistically convergent to $a$. For every $\eps>0$ and all $n\in \N\sm A$ we have 
\begin{align*}
&\sum_{k=1}^{\infty}\sum_{l=1}^{\infty}b_{nk}^{(i)}b_{nl}^{(i)}\chi_{D(s,\eps)}(k,l)\leq
\sum_{k=1}^{\infty}\sum_{l=1}^{\infty}b_{nk}^{(i)}b_{nl}^{(i)}\paren*{\chi_{D(s,a,\eps/2)}(k)+\chi_{D(s,a,\eps/2)}(l)} \\
&\leq 2M\sum_{k=1}^{\infty}b_{nk}^{(i)}\chi_{D(s,a,\eps/2)}(k) \to 0 \ \ \mathrm{along} \ I \ \mathrm{uniformly\ in} \ i\in S.
\end{align*}
\end{Proof}

The next two propositions are the analogues of \cite{connor4}*{Theorem 3}. Since their proofs parallel very much those of 
Proposition \ref{prop:strong stat} resp. \ref{prop:stat strong} they will be omitted. In the formulation of both propositions, 
we differ from our usual notation and allow $\F=(F_{kl}^{(i)})_{k,l\in \N,i\in S}$ to be a family in $\mathcal{M}\cup\mathcal{O}$ 
with index set $\N\times\N\times S$ instead of $\N\times S$.

\begin{proposition}\label{prop:suf pre cauchy}
Suppose that
\begin{equation*}
I\text{-}\lim\sum_{k=1}^{\infty}\sum_{l=1}^{\infty}b_{nk}^{(i)}b_{nl}^{(i)}F_{kl}^{(i)}(\abs*{s_k-s_l})=0 \ \ \mathrm{uniformly\ in} \ i\in S
\end{equation*}
and
\begin{equation*}
L(t):=\inf\set*{F_{kl}^{(i)}(t):k,l\in \N, i\in S}>0 \ \ \forall t>0.
\end{equation*}
Then $s$ is $\B^I$-statistically pre-Cauchy.
\end{proposition}

\begin{proposition}\label{prop:nec pre cauchy}
Suppose that $s$ is bounded and $\B^I$-statistically pre-Cauchy. If $\F$ is equicontinuous at $0$ and
\begin{equation*}
\exists A\in I\ M:=\sup\set*{\sum_{k=1}^{\infty}b_{nk}^{(i)}:n\in \N\sm A, i\in S}<\infty,
\end{equation*} 
as well as
\begin{equation*}
h(t):=\sup\set*{F_{kl}^{(i)}(t):k,l\in \N, i\in S}<\infty \ \ \forall t\geq0,
\end{equation*}
then we also have
\begin{equation*}
I\text{-}\lim\sum_{k=1}^{\infty}\sum_{l=1}^{\infty}b_{nk}^{(i)}b_{nl}^{(i)}F_{kl}^{(i)}(\abs*{s_k-s_l})=0 \ \ \mathrm{uniformly\ in} \ i\in S.
\end{equation*}
\end{proposition}

It was proved in \cite{connor4} that a statistically pre-Cauchy sequence $(s_n)_{n\in \N}$ which possesses a convergent 
subsequence $(s_{n_k})_{k\in \N}$ such that the set of indices $\set*{n_k:k\in \N}$ is ``large'' in the sense that
\begin{equation*} 
\liminf_{n\to \infty}\frac{1}{n}\abs*{\set*{n_k:k\in \N \ \mathrm{and} \  n_k\leq n}}>0
\end{equation*}
is statistically convergent. This result can be generalised in the following way.
\begin{theorem}\label{thm:pre cauchy subseq conv}
Suppose that $I\ssq J_{\B,I}$ and
\begin{equation*}
\sup\set*{\sum_{k=1}^{\infty}b_{nk}^{(i)}:n\in \N, i\in S}<\infty.
\end{equation*}
Let $a$ be any real or complex number. Let $s=(s_n)_{n\in \N}$ be a $\B^I$-statistically pre-Cauchy sequence 
and let $W\ssq \N$ be such that for every $\eps>0$ the set $\set*{n\in W: \abs*{s_n-a}\geq\eps}$ belongs to $I$ 
and furthermore
\begin{equation*}
w:=I\text{-}\liminf\inf_{i\in S}\sum_{k=1}^{\infty}b_{nk}^{(i)}\chi_W(k)>0.
\end{equation*}
Then $s$ is $\B^I$-statistically convergent to $a$.
\end{theorem}

\begin{Proof}
Take $\eps,\delta>0$ arbitrary. Then $V:=\set*{k\in W:\abs*{s_k-a}\geq\eps/2}\in I$, by assumption.
Put $A:=\set*{k\in W:\abs*{s_k-a}<\eps/2}$, $B:=\set*{k\in \N:\abs*{s_k-a}\geq\eps}$ and 
$C:=\set*{(k,l)\in \N^2:\abs*{s_k-s_l}\geq\eps/2}$. Then $A\times B\ssq C$.\par
Let us also fix $\tau\in (0,w)$ such that $\tau(w-\tau)^{-1}\leq\delta$. Since $s$ is $\B^I$-statistically
pre-Cauchy there is some $E\in I$ such that
\begin{equation*}
\set*{n\in \N:\sum_{k=1}^{\infty}\sum_{l=1}^{\infty}b_{nk}^{(i)}b_{nl}^{(i)}\chi_C(k,l)\geq\tau}\ssq E \ \ \forall i\in S.
\end{equation*}
But we have
\begin{equation*}
\sum_{k=1}^{\infty}\sum_{l=1}^{\infty}b_{nk}^{(i)}b_{nl}^{(i)}\chi_C(k,l)\geq
\paren*{\sum_{k=1}^{\infty}b_{nk}^{(i)}\chi_A(k)}\paren*{\sum_{l=1}^{\infty}b_{nl}^{(i)}\chi_B(l)}
\end{equation*}
and thus 
\begin{equation*}
\set*{n\in \N:\paren*{\sum_{k=1}^{\infty}b_{nk}^{(i)}\chi_A(k)}\paren*{\sum_{l=1}^{\infty}b_{nl}^{(i)}\chi_B(l)}\geq\tau}\ssq E \ \ \forall i\in S.
\end{equation*}
Since $V\in I\ssq J_{\B,I}$ it follows that
\begin{equation*}
I\text{-}\lim\sup_{i\in S}\sum_{k=1}^{\infty}b_{nk}^{(i)}\chi_V(k)=0.
\end{equation*}
Because of Lemma \ref{lemma:I limsup} this implies
\begin{align*}
&w=I\text{-}\liminf\inf_{i\in S}\sum_{k=1}^{\infty}b_{nk}^{(i)}(\chi_A(k)+\chi_V(k)) \\ 
&\leq I\text{-}\liminf\paren*{\inf_{i\in S}\sum_{k=1}^{\infty}b_{nk}^{(i)}\chi_A(k)+\sup_{i\in S}\sum_{k=1}^{\infty}b_{nk}^{(i)}\chi_V(k)} \\
&=I\text{-}\liminf\inf_{i\in S}\sum_{k=1}^{\infty}b_{nk}^{(i)}\chi_A(k)=:r.
\end{align*}
By \cite{demirci0}*{Theorem 2} we have
\begin{equation*}
F:=\set*{n\in \N:\inf_{i\in S}\sum_{k=1}^{\infty}b_{nk}^{(i)}\chi_A(k)<r-\tau}\in I.
\end{equation*}
If $n\in \N\sm(E\cup F)$ then $\sum_{k=1}^{\infty}b_{nk}^{(i)}\chi_B(k)<\tau(r-\tau)^{-1}\leq\tau(w-\tau)^{-1}\leq\delta$ for every $i\in S$.\par
Thus $E\cup F\in I$ with
\begin{equation*}
\set*{n\in \N:\sum_{k=1}^{\infty}b_{nk}^{(i)}\chi_B(k)\geq\delta}\ssq E\cup F \ \ \forall i\in S
\end{equation*}
and the proof is finished.
\end{Proof}

By \cite{connor4}*{Theorem 5} a bounded statistically pre-Cauchy sequence in $\R$ whose set of cluster points is
nowhere dense is statistically convergent. To obtain an analogous result in our setting, we introduce the following 
strengthening of the notion of $\B^I$-statistically pre-Cauchy sequences.
\begin{definition}\label{def:stat pre+ cauchy}
If each $B_i$ is non-negative, a sequence $s=(s_k)_{k\in \N}$ of real or complex numbers is called a $\B_+^I$-statistically
pre-Cauchy sequence if for every $\eps>0$
\begin{equation*}
I\text{-}\lim\sum_{k=1}^{\infty}\sum_{l=1}^{\infty}b_{nk}^{(i)}b_{nl}^{(j)}\chi_{D(s,\eps)}(k,l)=0 \ \ \mathrm{uniformly\ in} \ i,j\in S.
\end{equation*}
\end{definition}
For $\B_+^I$-statistically pre-Cauchy sequences, Lemma \ref{lemma:stat conv pre cauchy}, Proposition \ref{prop:suf pre cauchy} 
and Proposition \ref{prop:nec pre cauchy} hold accordingly (with the obvious modifications, one can even take a family 
$\F=(F_{kl}^{(i,j)})_{k,l\in \N,i,j\in S}$ in $\mathcal{M}\cup\mathcal{O}$ with index set $\N^2\times S^2$ in this case).\par
The next lemma generalises \cite{connor4}*{Lemma 4} while its proof follows the same lines.
\begin{lemma}\label{lemma:stat pre+ cauchy}
Let $I$ be an admissible ideal. Suppose that $\sum_{k=1}^{\infty}b_{nk}^{(i)}<\infty$ for all $n\in \N, i\in S$ and 
\begin{align}
&\exists A\in I\ M:=\sup\set*{\sum_{k=1}^{\infty}b_{nk}^{(i)}:n\in \N\sm A, i\in S}<\infty, \label{eq:5.1}\\
&I\text{-}\lim \sum_{k=1}^{\infty}b_{nk}^{(i)}=1 \ \ \mathrm{uniformly\ in} \ i\in S. \label{eq:5.2}
\end{align}
Let $\mathcal{W}$ be a basis for $\F(I)$ such that for every $\set*{n_1<n_2\dots n_k<n_{k+1}\dots}\in \mathcal{W}$
the following holds: 
\begin{equation}\label{eq:5.3}
\exists k_0\in \N \ \forall k\geq k_0 \ \inf_{i\in S}\sum_{l=1}^{\infty}\abs*{b_{n_kl}^{(i)}-b_{n_{k+1}l}^{(i)}}<\frac{1}{3}.
\end{equation}
Let $s=(s_n)_{n\in \N}$ be a $\B_+^I$-statistically pre-Cauchy sequence in $\R$ and $\alpha<\beta$ such that 
$H:=\set*{n\in \N:s_n\in (\alpha,\beta)}\in J_{\B,I}$.\par
Then $X:=\set*{n\in \N:s_n\leq\alpha}\in J_{\B,I}$ or $Y:=\set*{n\in \N:s_n\geq\beta}\in J_{\B,I}$.
\end{lemma}

\begin{Proof}
Let us put $t_n=s_n$ if $n\not\in H$ and $t_n=\alpha$ if $n\in H$. Since $H\in J_{\B,I}$, it is not difficult to 
see that $t=(t_n)_{n\in \N}$ is also $\B_+^I$-statistically pre-Cauchy. Put $P:=\set*{n\in \N:t_n\leq\alpha}$ and
$Q:=\set*{n\in \N:t_n\geq\beta}$. Then $X\ssq P\cup H$ and $Y\ssq Q\cup H$, thus it suffices to show $P\in J_{\B,I}$
or $Q\in J_{\B,I}$. Note also that $t_n\not\in (\alpha,\beta)$ for all $n\in \N$ and hence $Q=\N\sm P$.\par
For the sake of brevity, we define for $n\in \N$ and $i\in S$
\begin{equation*}
D_{ni}(K):=\sum_{k=1}^{\infty}b_{nk}^{(i)}\chi_A(k) \ \ \forall K\ssq \N.
\end{equation*}
We claim that
\begin{equation}\label{eq:5.4}
I\text{-}\lim D_{ni}(P)(1-D_{nj}(P))=0 \ \ \mathrm{uniformly\ in} \ i,j\in S.
\end{equation}
To see this, fix an arbitrary $\eps>0$ and note that $P\times Q\ssq D(t,\beta-\alpha)$. So, since $t$ is 
$\B_+^I$-statistically pre-Cauchy, there is some $E\in I$ such that
\begin{equation*}
\set*{n\in \N: D_{ni}(P)D_{nj}(Q)\geq \frac{\eps}{2}}\ssq E \ \ \forall i,j\in S.
\end{equation*}
By \eqref{eq:5.2} there exists $F\in I$ such that
\begin{equation*}
\set*{n\in \N:\abs*{D_{ni}(\N)-1}\geq\frac{\eps}{2M}}\ssq F \ \ \forall i\in S.
\end{equation*}
Because of \eqref{eq:5.1} and $D_{ni}(Q)=D_{ni}(\N)-D_{ni}(P)$ this easily implies
\begin{equation*}
\set*{n\in \N:\abs*{D_{ni}(P)(1-D_{nj}(P))}\geq\eps}\ssq E\cup F\cup A \ \ \forall i,j\in S,
\end{equation*}
proving our claim. In particular, we can find $C\in \mathcal{W}$ with
\begin{equation*}
\abs*{D_{ni}(P)(1-D_{nj}(P))}<\frac{1}{9} \ \ \forall n\in C, \forall i,j\in S.
\end{equation*}
Then for every $n\in C$ we must have
\begin{equation*}
\sup_{i\in S}D_{ni}(P)\leq\frac{1}{3} \ \ \mathrm{or} \ \ \inf_{j\in S} D_{nj}(P)\geq\frac{2}{3}.
\end{equation*}
Write $C=\set*{n_1<n_2\dots n_k<n_{k+1}\dots}$ and choose $k_0$ according to \eqref{eq:5.3}.
Suppose first that $\sup_{i\in S}D_{n_{k_0}i}(P)\leq 1/3$. Then the same must hold for every $k>k_0$,
for elsewhise we could find a minimal $k>k_0$ with $\inf_{i\in S}D_{{n_k}i}(P)\\ \geq 2/3$ which would imply
\begin{align*}
\sum_{l=1}^{\infty}\abs*{b_{n_kl}^{(i)}-b_{n_{k-1}l}^{(i)}}\geq D_{{n_k}i}(P)-D_{{n_{k-1}}i}(P)\geq
\frac{2}{3}-\frac{1}{3}=\frac{1}{3}
\end{align*}
for all $i\in S$, contradicting the choice of $k_0$.\par
So we have $D_{{n_k}i}(P)\leq 1/3$ for all $k\geq k_0$ and all $i\in S$. Now fix again an arbitray $\eps>0$.
By \eqref{eq:5.4} there is $G\in I$ such that
\begin{equation*}
\set*{n\in \N:\abs*{D_{ni}(P)(1-D_{nj}(P))}\geq \frac{2}{3}\eps}\ssq G \ \ \forall i,j\in S.
\end{equation*}
Since $I$ is admissible, $R:=G\cup (\N\sm \set*{n_k:k\geq k_0})$ is again an element of $I$ and we have
\begin{equation*}
\set*{n\in \N:D_{ni}(P)\geq\eps}\ssq R \ \ \forall i\in S.
\end{equation*}
Thus we have shown that $D_{ni}(P)$ converges along $I$ to zero uniformly in $i\in S$, which means exactly that
$P\in J_{\B,I}$.\par
In the second case, $\inf_{i\in S}D_{n_{k_0}i}(P)\geq 2/3$, one can show analogously that $Q\in J_{\B,I}$.
\end{Proof}

Note that if $I=I_f$ and $\inf_{i\in S}\sum_{l=1}^{\infty}\abs{b_{nl}^{(i)}-b_{n+1l}^{(i)}}<1/3$ for all but finitely 
many $n\in \N$, then we can take $\mathcal{W}=\set*{\set*{n\in \N:n\geq m}:m\in \N}$ and condition \eqref{eq:5.3} is 
satisfied. For the Ces\`aro-matrix $C$ we even have $\lim_{n\to \infty}\sum_{l=1}^{\infty}\abs*{c_{nl}-c_{n+1l}}=0$.\par
As in \cite{connor4}, we can now use the above lemma to obtain a sufficient condition for $\B^I$-statistical convergence.
\begin{theorem}\label{thm:stat pre+ cauchy}
Under the same general hypotheses as in the previous lemma, if $s=(s_n)_{n\in \N}$ is a $J_{\B,I}$-bounded 
$\B_+^I$-statistically pre-Cauchy sequence in $\R$ such that the set $Z$ of all $J_{\B,I}$-cluster points 
of $s$ is nowhere dense\footnote{Note that $Z$ is closed (cf.\,\cite{kostyrko1}*{Theorem 4.1(i)}), so 
``$Z$ nowhere dense'' just means that $Z$ has empty interior.} in $\R$, then $s$ is $\B^I$-statistically 
convergent.
\end{theorem}

\begin{Proof}
Suppose that $s$ is $J_{\B,I}$-bounded and $\B_+^I$-statistically pre-Cauchy but not $\B^I$-statistically 
convergent.\par
As mentioned before, the $J_{\B,I}$-boundedness assures that there is some $a\in Z$. Since $s$ is not
$\B^I$-statistically convergent there is an $\eps>0$ such that $\set*{n\in \N:s_n\leq a-\eps}\not\in J_{\B,I}$
or $\set*{n\in \N:s_n\geq a+\eps}\not\in J_{\B,I}$. Without loss of generality, assume the former.\par
As in \cite{connor4}, we will show that $(a-\eps,a)\ssq Z$. If not, there would be an open intervall 
$(\alpha,\beta)\ssq (a-\eps,a)$ such that $\set*{n\in \N:s_n\in (\alpha,\beta)}\in J_{\B,I}$.\par
It follows from Lemma \ref{lemma:stat pre+ cauchy} that $X=\set*{n\in \N:s_n\leq\alpha}\in J_{\B,I}$ or 
$Y:=\set*{n\in \N:s_n\geq\beta}\in J_{\B,I}$.\par
Since $X\supseteq \set*{n\in \N:s_n\leq a-\eps}\not\in J_{\B,I}$ we would have $Y\in J_{\B,I}$. But we can 
find $\delta>0$ with $\beta<a-\delta$ and because of $a\in Z$ the set $\set*{n\in \N: s_n>a-\delta}$
cannot belong to $J_{\B,I}$ where on the other hand it is contained in $Y$.\par
Thus $Z$ has non-empty interior and the proof is finished.
\end{Proof}

As an immediate consequence of Theorem \ref{thm:stat pre+ cauchy} we get the following corol\-lary.
\begin{corollary}\label{cor:stat pre+ cauchy}
Under the same general assumptions as in Lemma \ref{lemma:stat pre+ cauchy}, if $s$ is a $\B_+^I$-statistically 
pre-Cauchy sequence in $\R$ whose range is finite, then $s$ is $\B^I$-statistically convergent.
\end{corollary}

\section{A sup-limsup-theorem}\label{sec:sup limsup}
In this section we will present the generalisation of Simons' equality that was announced in the abstract, 
but first we need to recall some definitions: A boundary for a real Banach space $X$ is a subset $H$ of 
$B_{X^*}$\footnote{For every Banach space $Y$ we denote by $B_Y$ its closed unit ball and by $S_Y$ its unit sphere.}
such that for every $x\in X$ there is some $x^*\in H$ with $x^*(x)=\norm{x}$. By the Hahn-Banach-theorem, $S_{X^*}$ 
is always a boundary for $X$. It easily follows from the Krein-Milman-theorem that $\ex B_{X^*}$, the set of extreme 
points of $B_{X^*}$, is also a boundary for $X$.\par
A famous theorem due to Rainwater (cf.\,\cite{rainwater}) states that a bounded sequence in $X$ which is convergent 
to some $x\in X$ under every functional from $\ex B_{X^*}$ is weakly convergent to $x$.\par
Later Simons (cf.\,\cite{simons1} and \cite{simons2}) generalised this result to an arbitrary boundary $H$ by 
proving that for every bounded sequence $(x_n)_{n\in \N}$ in $X$ the equality
\begin{equation*}
\smashoperator{\sup_{x^*\in H}}\limsup x^*(x_n)=\smashoperator{\sup_{x^*\in B_{X^*}}}\limsup x^*(x_n),
\end{equation*}
which is nowadays known as Simons' equality, holds.\par
An easy separation argument shows that every boundary $H$ satisfies $B_{X^*}=\cco[w^*]{H}$, but $B_{X^*}=\cco{H}$
is not true in general (here $\co{A}$ denotes the convex hull, $\cl[w^*]{A}$ the weak*-closure and $\cl{A}$ the
norm-closure of $A\ssq X^*$).\par
In \cite{fonf} Fonf and Lindenstrauss introduced the following intermediate notion. Consider a convex weak*-compact 
subset $K$ of $X^*$ (where $X$ is a real or complex Banach space). A subset $H$ of $K$ is said to $(I)$-generate 
$K$ provided that whenever $H$ is written as a countable union $H=\bigcup_{m=1}^{\infty}H_m$ then
\begin{equation*}
\cco\paren*{\bigcup_{m=1}^{\infty}\cco[w^*]{H_m}}=K
\end{equation*}
or equivalently, whenever $H$ is written as a countable union $H=\bigcup_{m=1}^{\infty}H_m$ with $H_m\ssq H_{m+1}$
then
\begin{equation*}
\cl{\bigcup_{m=1}^{\infty}\cco[w^*]{H_m}}=K.
\end{equation*}
Clearly, $K=\cco{H}$ implies that $H$ $(I)$-generates $K$ which in turn implies $K=\cco[w^*]{H}$, but the converses are
not true in general as was shown in \cite{fonf}. It was also proved in \cite{fonf} that, for a real Banach space, every 
boundary of $K$ $(I)$-generates $K$.\footnote{The set $H$ is called a boundary of $K$ if $\max\set*{x^*(x):x^*\in H}=
\sup\set*{x^*(x):x^*\in K}$ for every $x\in X$. In this terminology, $H$ is a boundary for $X$ \ifif it is a boundary 
of $B_{X^*}$.}\par
Nygaard proved in \cite{nygaard} that Rainwater's theorem holds true for every $(I)$-generating subset of $B_{X^*}$ and
the authors of \cite{cascales} showed that Simons' equality is equivalent to the $(I)$-generation property
(cf.\,\cite{cascales}*{Theorem 2.2}, see also \cite{kalenda1}*{Lemma 2.1 and Remark 2.2}).\par
In \cite{hardtke} the author investigated the possibility to generalise the Rainwater-Simons-convergence theorem for
$(I)$-generating sets to some generalised convergence methods such as strong $A$-$\mathbf{p}$-summability and almost 
convergence by proving a general Simons-like inequality for $(I)$-generating sets (cf.\,\cite{hardtke}*{Theorem 3.1}).
We will continue this work here, using similiar arguments as in \cite{hardtke} to generalise Simons' equality to
the $J_{\B,I}\text{-}\limsup$ for the case that $\F(I)$ has a countable base and obtain some related convergence results.\par
First we need the following lemma, whose proof is---once more---analo\-gous to those of the Propositions \ref{prop:strong stat} 
and \ref{prop:stat strong}. Therefore, the details will be skipped.
\begin{lemma}\label{lemma:aux 2}
Let each $B_i$ be non-negative. Define $f:\R \rightarrow [0,\infty)$ by $f(t)=t$ for $t\geq 0$ and $f(t)=0$ for $t<0$. Put
$A(s,a,\eps):=\set*{k\in \N:s_k>a+\eps}$ for every $\eps>0$. Then
\begin{align*}
&I\text{-}\lim\sum_{k=1}^{\infty}b_{nk}^{(i)}f(s_k-a)=0 \ \ \mathrm{uniformly\ in} \ i\in S \\
&\Rightarrow \ A(s,a,\eps)\in J_{\B,I} \ \ \forall \eps>0
\end{align*}
and the converse is true if the sequence $s$ is bounded and
\begin{equation*}
\sup\set*{\sum_{k=1}^{\infty}b_{nk}^{(i)}:n\in \N\sm A, i\in S}<\infty
\end{equation*}
for some $A\in I$.
\end{lemma}

Now for the generalisation of Simons' equality.
\begin{theorem}\label{thm:simons eq}
Let $X$ be a real Banach space, $K\ssq X^*$ a convex weak*-compact subset and $H\ssq K$ an $(I)$-generating set for $K$. 
Let the ideal $I$ be such that the filter $\F(I)$ has a countable base. assume that Each $B_i$ is non-neagtive and that 
there exists an $A\in I$ such that
\begin{equation*}
M:=\sup\set*{\sum_{k=1}^{\infty}b_{nk}^{(i)}:n\in \N\sm A, i\in S}<\infty.
\end{equation*}
Let $(x_n)_{n\in \N}$ be a bounded sequence in $X$. Then the equality
\begin{equation*}
\smashoperator{\sup_{x^*\in H}}J_{\B,I}\text{-}\limsup x^*(x_n)=\smashoperator{\sup_{x^*\in K}}J_{\B,I}\text{-}\limsup x^*(x_n)
\end{equation*}
holds.
\end{theorem}

\begin{Proof}
Denote the left-hand supremum by $c$, the right-hand supremum by $d$. We only have to show $d\leq c$. Let $R=\sup_{n\in \N}\norm{x_n}$.
Let $(C_n)_{n\in \N}$ be a countable base for $\F(I)$. Without loss of generality we may assume $C_{n+1}\ssq C_n$ for all $n$. Take
$x^*\in K$ and $\eps>0$ arbitrary and put
\begin{align*}
&E_m=\set*{y^*\in K:\sum_{k=1}^{\infty}b_{nk}^{(i)}f(y^*(x_k)-c)\leq\eps \ \forall i\in S,n\in C_m} \\ 
&\mathrm{and} \ H_m=E_m\cap H \ \ \forall m\in \N,
\end{align*}
where $f$ is as in the previous lemma. Then $H_m\ssq H_{m+1}$ for every $m\in \N$. It follows from \cite{demirci0}*{Theorem 1} 
that $\set*{n\in \N:y^*(x_n)>c+\delta}\in J_{\B,I}$ for every $\delta>0$. Together with the previous lemma this easily implies 
$\bigcup_{m=1}^{\infty}H_m=H$.\par
Since $H$ $(I)$-generates $K$ we get that
\begin{equation*}
K=\cl{\bigcup_{m=1}^{\infty}\cco[w^*]{H_m}}.
\end{equation*}
Thus we can find $m\in \N$ and $y^*\in \cco[w^*]{H_m}$ with $\norm{x^*-y^*}\leq\eps$. It is easily checked that $E_m$ is convex and weak*-closed,
hence $y^*\in E_m$. But for every $k\in \N$
\begin{align*}
&f(x^*(x_k)-c)\leq f(x^*(x_k)-y^*(x_k))+f(y^*(x_k)-c) \\
&\leq \norm{x^*-y^*}\norm{x_k}+f(y^*(x_k)-c)\leq R\eps+f(y^*(x_k)-c).
\end{align*}
It follows that
\begin{align*}
\sum_{k=1}^{\infty}b_{nk}^{(i)}f(x^*(x_k)-c)\leq MR\eps+\sum_{k=1}^{\infty}b_{nk}^{(i)}f(y^*(x_k)-c)\leq\eps(MR+1)
\end{align*}
for every $i\in S$ and every $n\in C_m\cap(\N\sm A)$. Since $C_m\cap(\N\sm A)\in \F(I)$ and $\eps>0$ was arbitrary we 
conclude with Lemma \ref{lemma:aux 2} that $\set*{n\in \N:x^*(x_n)>c+\delta}\in J_{\B,I}$ for every $\delta>0$, whence 
$J_{\B,I}\text{-}\limsup x^*(x_n)\leq c$.
\end{Proof}

As a corollary, we get the following convergence result.
\begin{corollary}\label{cor:simons conv}
Under the same hypotheses as in Theorem \ref{thm:simons eq} with $K=B_{X^*}$, if $x\in X$ 
is such that $(x^*(x_n))_{n\in \N}$ is $\B^I$-statistically convergent to $x^*(x)$ for every 
$x^*\in H$ then the same holds true for every $x^*\in X^*$, i.\,e., $(x_n)_{n\in \N}$
is ``weakly $\B^I$-statistically convergent to $x$''.\par
Moreover, for every family $\F=(F_k^{(i)})_{k\in\N,i\in S}$ in $\mathcal{M}\cup\mathcal{O}$ 
which is equicontinuous at $0$ and satisfies
\begin{equation*}
\inf\set*{F_k^{(i)}(t):k\in \N, i\in S}>0 \ \ \forall t>0
\end{equation*}
and
\begin{equation*}
\sup\set*{F_k^{(i)}(t):k\in \N, i\in S}<\infty \ \ \forall t\geq0,
\end{equation*}
$(x^*(x_n))_{n\in \N}$ is strongly $\B^I$-summable to $x^*(x)$ with respect to $\F$ 
for every $x^*\in X^*$ whenever this statement holds for every $x^*\in H$.
\end{corollary}

\begin{Proof}
The first statement follows directly from Theorem \ref{thm:simons eq} and the second follows 
from the first one via the Propositions \ref{prop:strong stat} and \ref{prop:stat strong}.
\end{Proof}

It is clear that this convergence result carries over to complex Banach spaces (note that if $X$ 
is a complex Banach space and $H$ $(I)$-generates $B_{X^*}$ then $\set*{\RE x^*:x^*\in H}$ 
$(I)$-generates $\set*{\RE x^*:x^*\in B_{X^*}}$, the unit ball of the underlying real space).\par
In particular, if we take each $B_i$ to be the infinite unit matrix, we get that for every ideal 
$I$ such that $\F(I)$ has a countable base, $I\text{-}\lim x^*(x_n)=x^*(x)$ for every $x^*\in X^*$ 
whenever this is true for every $x^*$ in an $(I)$-generating subset of $B_{X^*}$ (in particular,
in a boundary for $X$). We can also prove an analogous convergence result for $\B^I$-summability.
\begin{proposition}\label{prop:simons BI-sum}
Let $X$ be a real or complex Banach space and $H\ssq B_{X^*}$ an $(I)$-generating set for $B_{X^*}$.
Suppose that $\F(I)$ has a countable base, $\sum_{k=1}^{\infty}\abs{b_{nk}^{(i)}}<\infty$ for all 
$n\in \N, i\in S$ and moreover
\begin{equation*}
M:=\sup\set*{\sum_{k=1}^{\infty}\abs*{b_{nk}^{(i)}}:n\in \N\sm A, i\in S}<\infty
\end{equation*}
for some $A\in I$.\par
Let $(x_n)_{n\in \N}$ be a bounded sequence in $X$ and $x\in X$ such that $(x^*(x_n))_{n\in \N}$
is $\B^I$-summable to $x^*(x)$ for every $x^*\in H$. Then the same is true for every $x^*\in X^*$.
\end{proposition}

\begin{Proof}
Let $(C_n)_{n\in \N}$ be a decreasing countable basis for $\F(I)$. Let $R\geq\sup_{n\in \N}\norm{x_n}$
and $R\geq\norm{x}$. Take any $x^*\in B_{X^*}$ and fix an arbitrary $\eps>0$. Define
\begin{align*}
&E_m:=\set*{y^*\in B_{X^*}:\sup_{i\in S}\abs*{\sum_{k=1}^{\infty}b_{nk}^{(i)}y^*(x_k)-y^*(x)}\leq\eps \ \forall n\in C_m} \\
&\mathrm{and} \ H_m:=E_m\cap H \ \ \forall m\in \N.
\end{align*}
Then $H_m\nearrow H$ and since $H$ $(I)$-generates $B_{X^*}$ we can find $m\in \N$ and $y^*\in \cco[w^*]{H_m}$
such that $\norm{x^*-y^*}\leq\eps$.\par
It is not too hard to see that $E_m$ is convex and weak*-closed and thus $y^*\in E_m$. Consequently, for all 
$i\in S$ and $n\in C_m\cap(\N\sm A)$ we have
\begin{align*}
&\abs*{\sum_{k=1}^{\infty}b_{nk}^{(i)}x^*(x_k)-x^*(x)}\leq\abs*{\sum_{k=1}^{\infty}b_{nk}^{(i)}(x^*(x_k)-y^*(x_k))} \\
&+\abs*{\sum_{k=1}^{\infty}b_{nk}^{(i)}y^*(x_k)-y^*(x)}+\abs*{y^*(x)-x^*(x)} \\
&\leq M\norm{x^*-y^*}R+\eps+\norm{x^*-y^*}R\leq\eps(R(M+1)+1).
\end{align*}
Since $C_m\cap(\N\sm A)\in \F(I)$ and $\eps>0$ was arbitrary we are done.
\end{Proof}

The next result concerning $\B^I$-statistically pre-Cauchy sequences is a generalisation of \cite{hardtke}*{Corollary 3.5}.
Using Proposition \ref{prop:suf pre cauchy} and Proposition \ref{prop:nec pre cauchy} with $F_{kl}^{(i)}=\id_{[0,\infty)}$ 
for all $k,l\in \N$ and $i\in S$ its proof can be carried out analogously to that of Proposition \ref{prop:simons BI-sum}. 
The details will be omitted. 
\begin{proposition}\label{prop:simons pre cauchy}
Let $X$ be a real or complex Banach space and $H\ssq B_{X^*}$ an $(I)$-generating set for $B_{X^*}$.
Suppose that $\F(I)$ has a countable base, that each $B_i$ is non-negative and that there is some $A\in I$
such that
\begin{equation*}
\sup\set*{\sum_{k=1}^{\infty}b_{nk}^{(i)}:n\in \N\sm A, i\in S}<\infty.
\end{equation*}
Let $(x_n)_{n\in \N}$ be a bounded sequence in $X$ such that $(x^*(x_n))_{n\in \N}$ 
is $\B^I$-statisti\-cally pre-Cauchy resp. $\B_+^I$-statistically pre-Cauchy for every 
$x^*\in H$. Then the same is true for every $x^*\in X^*$.
\end{proposition}

Finally, let us give characterisations of weak-compactness and reflexivity that generalise 
\cite{hardtke}*{Corollaries 3.7 and 3.8}.
\begin{corollary}\label{cor:char weak comp}
Let $M$ be a bounded subset of the Banach space $X$ and $B$ an $(I)$-generating set for $B_{X^*}$. Then $M$ is 
weakly relatively compact if (and only if) for every sequence $(x_n)_{n\in \N}$ in $M$ there is an element 
$x\in X$, an ideal $I$ on $\N$ such that $\F(I)$ admits a countable base and a non-negative matrix 
$A=(a_{nk})_{n,k\geq 1}$ such that
\begin{align}
&\exists C\in I\ \sup_{n\in \N\sm C}\sum_{k=1}^{\infty}a_{nk}<\infty, \label{eq:6.1}\\
&I\text{-}\lim a_{nk}=0 \ \ \forall k\in \N \label{eq:6.2}
\end{align}
and $(x^*(x_n))_{n\in \N}$ is $A^I$-statistically convergent to $x^*(x)$ for every $x^*\in B$.
\end{corollary}

\begin{Proof}
Let $(x_n)_{n\in \N}$ be an arbitrary sequence in $M$ and fix $x$, $I$ and $A$ as above.
By Corollary \ref{cor:simons conv} $(x^*(x_n))_{n\in \N}$ is $A^I$-statistically convergent 
to $x^*(x)$ for every $x^*\in X^*$. Thus, given finitely many functionals $x_1^*,\dots,x_m^*\in X^*$,
the sequence $(\sum_{j=1}^m\abs{x_j^*(x_n-x)})_{n\in \N}$ is $A^I$-statistically convergent to zero.
Hence for any $\eps>0$ the set $D_{\eps}=\set*{n\in \N:\sum_{j=1}^m\abs{x_j^*(x_n-x)}<\eps}$ does not 
belong to $J_{A,I}$.\par
By \eqref{eq:6.2}, $J_{A,I}$ is admissible, therefore $D_{\eps}$ must be infinite for every $\eps>0$, 
which shows that $x$ is a weak-cluster point of $(x_n)_ {n\in \N}$.\par
So $M$ is weakly relatively countably compact and by the Eberlein-Shmulyan theorem, it must be also 
weakly relatively compact.
\end{Proof}

\begin{corollary}\label{cor:char reflex}
If $B_X$ is an $(I)$-generating set for $B_{X^{**}}$\footnote{We consider $X$ canonically embedded 
into its bidual.}, then $X$ is reflexive if (and only if) for every sequence $(x_n^*)_{n\in \N}$ in 
$B_{X^*}$ there is a functional $x^*\in X^*$, an ideal $I$ on $\N$ such that $\F(I)$ admits a 
countable base and a non-negative matrix $A$ such that \eqref{eq:6.1} and \eqref{eq:6.2} are 
satisfied and $(x_n^*(x))_{n\in \N}$ is $A^I$-statistically convergent to $x^*(x)$ for every $x\in X$.
\end{corollary}

\begin{Proof}
By the previous corollary, $B_{X^*}$ is weakly compact, thus $X^*$ and hence also $X$ is reflexive.
\end{Proof}

\address
\email

\end{document}